\documentclass{amsart}
\usepackage{amssymb}
\usepackage{mathrsfs}
\usepackage{bbm}
\usepackage{oldgerm}
\usepackage[english]{babel}
\usepackage[T1]{fontenc}
\usepackage[latin1]{inputenc}
\usepackage[all]{xy}
\usepackage{hyperref}
\usepackage{upref}
\usepackage{xcolor}
\usepackage{tikz-cd}

\hypersetup{
 colorlinks,
 linkcolor={red!50!black},
 citecolor={blue!50!black},
 urlcolor={blue!80!black}
}

\newtheorem{teo}[subsection]{Theorem}
\newtheorem{prop}[subsection]{Proposition}
\newtheorem{cor}[subsection]{Corollary}
\newtheorem{lem}[subsection]{Lemma}

\theoremstyle{definition}

\newtheorem{rema}[subsection]{Remark}

\numberwithin{equation}{subsection}

\newcommand{\Aut}{\mathrm{Aut}}
\newcommand{\Art}{\mathrm{Art}}

\newcommand{\colim}{\mathrm{colim}}

\newcommand{\F}{\mathcal{F}}
\newcommand{\Frob}{\mathrm{Frob}}

\newcommand{\G}{\mathcal{G}}

\newcommand{\Gal}{\mathrm{Gal}}
\newcommand{\Hom}{\mathrm{Hom}}
\newcommand{\Ind}{\mathrm{Ind}}
\newcommand{\id}{\mathrm{id}}

\newcommand{\Lc}{\mathcal{L}}
\newcommand{\Mac}{\mathcal{M}}

\newcommand{\Ow}{\mathcal{O}}

\newcommand{\rk}{\mathrm{rk}}
\newcommand{\Sh}{\mathrm{Sh}}

\newcommand{\Spec}{\mathrm{Spec}}

\newcommand{\Tr}{\mathrm{Tr}}

\newcommand{\Z}{\mathbb{Z}}

\title{Compatible systems of $\ell$-adic sheaves}
\author{Quentin Guignard}
\address{Max Planck Institute for Mathematics, Vivatsgasse 7, 53111 Bonn, Germany}
\email{guignard@mpim-bonn.mpg.de}

\begin{document}

\begin{abstract} We introduce a notion of compatibility for families $(\F_{\ell})_{\ell}$ of bounded constructible $\ell$-adic complexes of \'etale sheaves on schemes. For schemes of finite type over a field, this notion is preserved by the usual six functors. We prove that the compatibility of a family is preserved by the nearby cycles functor and by the linearized $\varepsilon$-factors from \cite{G192}. We establish independence of $\ell$ for the characteristic cycles and characteristic $\varepsilon$-cycles of compatible families.
 
\end{abstract}

\maketitle
\tableofcontents

\section{Introduction \label{intro}}

\subsection{\label{intro1}} In \cite[5.4]{Ill10}, Illusie formulates the problem of associating to any reasonable scheme $X$ a triangulated category
$$
M(X,\mathbb{Q}) \subseteq \prod_{\substack{\ell}} D_c^b \left(X\left[ \frac{1}{\ell} \right], \mathbb{Q}_{\ell} \right),
$$
of families of \'etale $\ell$-adic complexes on $X$, satisfying the six functors formalism, i.e. stability by $f^*,f^!,Rf_*,Rf_!, RHom, \otimes^L$, and satisfying properties of uniform constructibility similar to those proved in \cite[2.5-2.7]{Ill10}. For example, for any family $(\F_{\ell})_{\ell}$ in $M(X,\mathbb{Q})$, we would like to ask for the existence of a dense open subset $U$ of $X$ such that for each $\ell$, the cohomology sheaves of $\F_{\ell}$ are lisse on $U\left[ \frac{1}{\ell} \right]$.

\subsection{\label{intro3}} More generally, if $E$ is a number field, one can similarly ask for triangulated categories $M(X,E)$ in $\prod_{\substack{\lambda}} D_c^b \left(X\left[ \ell_{\lambda}^{-1} \right], E_{\lambda} \right)$ satisfying the same stabilities and uniform constructibility properties, where $\lambda$ runs over all non archimedean places of $E$ and where $\ell_{\lambda}$ is the residual characteristic of a place $\lambda$.

\subsection{\label{intro2}} We consider in this text the smallest collection of triangulated categories $M(X,E)$ as above, satisfying the following three conditions:
\begin{enumerate}
\item for any object $\F$ of $M(X,E)$, the Tate twist $\F(1)$ is in $M(X,E)$ as well;
\item for any object $\F$ in $D_c^b(X,E)$, the family $(\F \otimes_{E}  E_{\lambda})_{\lambda}$  is in $M(X,E)$;
\item for any separated morphism of finite presentation $f : X \rightarrow Y$, the functor $Rf_!$ sends $M(X,E)$ to $M(Y,E)$.
\end{enumerate}
We refer to \ref{3.24}, \ref{3.19}, \ref{6.1} for precise definitions. We note that these categories contain most complexes of interest in the study of exponential sums over finite fields. Indeed, sheaves of K\"ummer or Artin-Schreier type are already in $D_c^b(X,E)$ for a suitable number field $E$, hence belong to $M(X,E)$. 

\subsection{\label{intro5}} We prove that these triangulated categories $M(X,E)$ have the following properties:
\begin{enumerate}
\item uniform constructibility holds for families in $M(X,E)$ (cf. \ref{prop3.23.3});
\item stability by $Rf_!,f^*,\otimes$ (cf. \ref{prop3.23.2}, \ref{prop3.19.10});
\item stability by the vanishing and nearby cycles functors (cf. \ref{teo3.3.1}, \ref{3.20});
\item for schemes of finite type over a field, stability by $Rf_*,f^!, RHom$ (cf. \ref{prop6.3.1}, \ref{cor6.2.2}, \ref{cor6.2.3}).
\end{enumerate}
We refer to the indicated sections for precise statements. We more generally construct triangulated categories $M(X,Q,E)$ of $Q$-equivariant $E$-linear coefficients, for a profinite group $Q$ acting admissibly on $X$. This additional generality is useful in certain applications, but is also needed to state and prove the stability by the nearby cycles functor. We should note that the tools used to prove the latter stability, including the equivariant version of de Jong's semistable alteration theorem, were developped in Vidal's seminal article \cite{Vi} in order to generalize the results from \cite{Ill81}.

\subsection{\label{intro4}} For a scheme $X$ of finite type over $\mathbb{Z}$, families $(\F_{\lambda})_{\lambda}$ in $M(X,E)$ are $E$-compatible in the classical sense, i.e. for any closed point $x$ of $X$, the local $L$-function
$$
\det(1 - T \Frob_x | \F_{\lambda, \overline{x}}),
$$
is a rational function with coefficients in $E$, which is independent of $\lambda$ (cf. \ref{invar}). For schemes of finite type over a finite field, this notion of compatibility is already known to be preserved by the six functors $f^*,f^!,Rf_*,Rf_!, RHom, \otimes^L$, by results of Gabber, cf. Fujiwara's account of his proofs in \cite{Fu02}. Similar results in the context of schemes of finite type over a local field have been obtained by Zheng \cite{Zh09}. The corresponding stabilities of $M(X,E)$ can be considered as a refinement of these results: $E$-compatibility is a property of virtual sheaves, while belonging to $M(X,E)$ is a property of the actual complexes.

Similarly, the aforementioned stabilities of $M(X,E)$, together with the $E$-compatibility of families in $M(X,E)$ or $M(X,Q,E)$, form a refinement of various $\ell$-independence results, such as those of Ochiai \cite[Th. B]{Oc99}, Vidal \cite[Prop. 4.2]{Vi}, or Kato \cite[Prop. 5.2]{Ka19}.

\subsection{\label{intro6}} Let $X$ be a smooth scheme of pure dimension $d$ over a perfect field $k$ and let $\lambda$ be a non archimedean place of $E$ with residual characteristic invertible in $k$. Given an object $\F_{\lambda}$ of $D_{c}^b \left( X , \Ow_{E_{\lambda}} \right)$, one can consider Beilinson's singular support $\mathrm{SS}(\F_{\lambda})$ \cite{Bei}, and Saito's characteristic cycle $\mathrm{CC}(\F_{\lambda})$ \cite{Sai}. The singular support $\mathrm{SS}(\F_{\lambda})$ is a closed conical subset of the cotangent scheme $T^*X$ with irreducible components of dimension $d$, and $\mathrm{CC}(\F_{\lambda})$ is a $d$-cycle on $T^*X$ supported on $\mathrm{SS}(\F_{\lambda})$, characterized by the following property: for any commutative diagram 
\begin{center}
 \begin{tikzpicture}[scale=1]
\node (B) at (1,1) {$u$};
\node (A) at (2,0) {$S$};
\node (C) at (4,0) {$s$};
\node (D) at (2,1) {$U$};
\node (E) at (4,1) {$X$};

\path[->,font=\scriptsize]
(B) edge  (D)
(A) edge  (C)
(D) edge node[left]{$f$} (A)
(E) edge (C)
(D) edge node[above]{$j$} (E);
\end{tikzpicture} 
\end{center}
where $j$ is \'etale, $u$ is a closed point of $U$ and $S$ is a smooth curve over $s$, such that $df^{-1}(\mathrm{SS}(\F_{\lambda})_{|U \setminus \{ u\}})$ is contained in the $0$-section of $T^*S$, we have the Milnor formula
$$
- a_0( S_{(f(\overline{u}))}, R\Phi_{f}(\F_{\lambda})_{\overline{u}} ) = (j^*\mathrm{CC}(\F_{\lambda}),df)_{T^*U,u}
$$
where $a_0$ is the Artin conductor, or total dimension. The characteristic cycle $\mathrm{CC}(\F_{\lambda})$ depends only on the isomorphism class of $\F_{\lambda}$ in $D_{c}^b \left( X , E_{\lambda} \right)$. We prove that the characteristic cycle is well defined for families in $M(X,E)$.

\begin{teo}(Cor. \ref{invarcharccycle})\label{invarcharcycle2} Let $E$ be a number field. Let $X$ be a smooth scheme over a perfect field $k$, and let $\F = (\F_{\lambda})_{\lambda}$ be an object of $M(X,E)$. Then the characteristic cycle $\mathrm{CC}(\F_{\lambda})$ does not depend on the choice of a non archimedean place of $E$ with residual characteristic invertible in $k$.
\end{teo}

We provide a self-contained proof of this assertion, but this result can also be recovered from the $\ell$-independence of characteristic cycles proved in \cite{SY} and \cite{Ka19}.

\subsection{\label{intro7}} Let $k$ be a perfect field of positive characteristic $p$, and let $E$ be a number field containing a non trivial $p$-th root of unity. Let $X$ be a smooth $k$-scheme of pure dimension $d$ and let $\lambda$ be a non archimedean place of $E$ which does not divide $p$. Given an object $\F_{\lambda}$ of $D_{c}^b \left( X , \Ow_{E_{\lambda}} \right)$, one can consider Takeuchi's characteristic $\varepsilon$-cycle $\mathcal{E}(\F_{\lambda})$ \cite[Th. 4.9]{Ta}. The latter is a $d$-dimensional cycle on $T^*X$, supported on the singular support $\mathrm{SS}(\F_{\lambda})$, with coefficients in the abelian group 
$$
\mathbb{Q} \otimes_{\mathbb{Z}} \Hom_{\mathrm{cont}}(\Gal(\overline{k}/k) \rightarrow \Ow_{E_{\lambda}}^{\times}/\mu_{E_{\lambda}}),
$$
where $\mu_{E_{\lambda}}$ is the group of roots of unity in $E_{\lambda}$, and is characterized by the following property: for any commutative diagram
\begin{center}
 \begin{tikzpicture}[scale=1]
\node (B) at (1,1) {$u$};
\node (A) at (2,0) {$S$};
\node (C) at (4,0) {$s$};
\node (D) at (2,1) {$U$};
\node (E) at (4,1) {$X$};

\path[->,font=\scriptsize]
(B) edge  (D)
(A) edge  (C)
(D) edge node[left]{$f$} (A)
(E) edge (C)
(D) edge node[above]{$j$} (E);
\end{tikzpicture} 
\end{center}
as in \ref{intro6}, we have
$$
\varepsilon_0( S_{(f(u))}, R\Phi_{f}(\F_{\lambda})_{\overline{u}} )^{-1} \circ \mathrm{Ver}_{u/k} = (j^*\mathrm{CC}(\F_{\lambda}),df)_{T^*U,u}^{[k(u):k]},
$$
where $\varepsilon_0( S_{(f(u))},-)$ is the geometric $\varepsilon$-factor from \cite{G19}, considered as a continuous homomorphism from $\Gal(\overline{k}/k)$ to $\Ow_{E_{\lambda}}/\mu_{E_{\lambda}}$,  and where $\mathrm{Ver}_{u/k} :   \Gal(\overline{k}/k)^{\mathrm{ab}} \rightarrow \Gal(\overline{u}/u)^{\mathrm{ab}}$ is the transfer homomorphism. The cycle $\mathcal{E}(\F_{\lambda})$ depends only on the isomorphism class of $\F_{\lambda}$ in $D_{c}^b \left( X , E_{\lambda} \right)$.

\begin{teo}(Cor. \ref{invarepscycle})\label{invarepscycle2} Let $k$ be a perfect field of positive characteristic $p$ with separable closure $\overline{k}$, and let $E$ be a number field containing a non trivial $p$-th root of unity. Let $X$ be a smooth $k$-scheme and let $\F = (\F_{\lambda})_{\lambda}$ be an object of $M(X,E)$. Then there exists a unique cycle $\mathcal{E}(\F)$ on $T^* X$, with coefficients in 
$$
\mathbb{Q} \otimes_{\mathbb{Z}} \Hom_{\mathrm{cont}} \left(\Gal(\overline{k}/k) \rightarrow \widehat{\Ow_{E}\left[ p^{-1} \right]}^{\times}/\mu_{E} \right),
$$
where $\widehat{\Ow_{E}\left[ p^{-1} \right]}^{\times}$ is the profinite completion of the finitely generated abelian group $\Ow_{E}\left[ p^{-1} \right]^{\times}$, such that for any non archimedean place $\lambda$ of $E$ not dividing $p$, the $\varepsilon$-cycle $\mathcal{E}(\F_{\lambda})$ is obtained from $\mathcal{E}(\F)$ through the natural homomorphism
$$
 \widehat{\Ow_{E}\left[ p^{-1} \right]}^{\times}/\mu_{E} \rightarrow  \Ow_{E_{\lambda}}^{\times}/\mu_{E_{\lambda}} .
$$
\end{teo}

If the base field $k$ is the perfection of a finitely generated extension of $\mathbb{F}_p$, then we will prove that the coefficients of the cycle $\mathcal{E}(\F)$ belong to the smaller abelian group
$$
\mathbb{Q} \otimes_{\mathbb{Z}} \Hom_{\mathrm{cont}} \left(W(\overline{k}/k) \rightarrow \Ow_{E}\left[ p^{-1} \right]^{\times}/\mu_{E} \right),
$$
where $W(\overline{k}/k)$ is the Weil group of $k$. In particular, if $k$ is finite, then these coefficients can be interpreted as elements of the finite dimensional $\mathbb{Q}$-vector space $\mathbb{Q} \otimes_{\mathbb{Z}} \Ow_{E}\left[ p^{-1} \right]^{\times}$.

\subsection{\label{intro9}} Let $S$ be a henselian trait of positive characteristic $p$, with generic point $\eta$, perfect closed point $s$ and uniformizer $\pi$, and let $\psi : \mathbb{F}_p \rightarrow E^{\times}$ be a non trivial homomorphism. In \cite{G192}, we introduced triangulated functors
$$
\mathrm{Art}_{\pi,\psi} R \Phi_{X/S} : D_c^b(X, E_{\lambda}) \rightarrow D_c^b(X_{s}, E_{\lambda}),
$$
for schemes of finite type over $S$, which are compatible with proper pushforward and smooth pullbacks, and which linearize Artin conductors and local $\varepsilon$-factors for $X = S$, in the following sense: for any $\F$ in $D_c^b(S, E_{\lambda})$, we have
\begin{align*}
\rk(\mathrm{Art}_{\pi,\psi} R \Phi_{S/S} \F)) &= a(S,\F), \\
\det(\mathrm{Art}_{\pi,\psi} R \Phi_{S/S} \F) &= \varepsilon(S,\F,d\pi),
\end{align*}
where $a$ is the Artin conductor and $\varepsilon$ is the geometric $\varepsilon$-factor from \cite{G19}.

\begin{teo}(Prop. \ref{propvareps})\label{compati} Let $S$ be a henselian trait of positive characteristic $p$, with perfect closed point $s$, and uniformizer $\pi$. Let $E$ be a number field and let $\psi : \mathbb{F}_p \rightarrow E^{\times}$ be a non trivial homomorphism. Then for any $S$-scheme $X$ of finite type, the functor $\mathrm{Art}_{\pi,\psi} R \Phi_{X/S}$ sends $M(X,E)$ to $M(X_s,E)$.
\end{teo}

This result is a key step in the proof of Theorem \ref{invarepscycle2}.

\subsection{\label{intro8}} This article is organized as follows. In Section \ref{compat}, we introduce the triangulated categories $M^{\mathrm{eff}}(X,\Ow_E)$, or rather their equivariant versions $M^{\mathrm{eff}}(X,Q,\Ow_E)$, from which $M(X,Q,E)$ is deduced by localization and inversion of the Tate twist, and we prove their stability by the functors $Rf_!, f^*, \otimes_{\Ow_E}^L$, cf. \ref{prop3.23.2}, \ref{prop3.19.10}. We then prove in Section \ref{nearby} that the categories $M(X,Q,E)$ are preserved by the nearby and vanishing cycles functors, by using Gabber-Vidal's equivariant version of de Jong's semistable alteration theorem. The latter result is used in Section \ref{full} in order to obtain the full six functors formalism for $M(X,Q,E)$ for schemes of finite type over a field. We then prove Theorem \ref{compati}, or rather an equivariant version thereof, in Section \ref{vareps}. Finally, we establish $E$-compatibility for families in $M(X,Q,E)$ and prove Theorems \ref{invarcharcycle2} and \ref{invarepscycle2} in Section \ref{indep}.

\subsection{Acknowledgements} The author prepared this article while benefiting from a postdoctoral fellowship from the Max Planck Institute for Mathematics.

\section{Compatible systems \label{compat}}

For any scheme $X$ and any integer $n \geq 1$, we denote by $X[\frac{1}{n}]$ the fiber product $X \times_{\Spec(\mathbb{Z})} \Spec(\mathbb{Z}[\frac{1}{n}])$. Throughout this section, we fix a number field $E$, with ring of integers $\Ow_E$.

\subsection{\label{3.4}} A right action of a profinite group $Q$ on a scheme $X$ is said to be \textit{admissible} if for any point $x$ of $X$, there exists an affine open subset $U$ of $X$ containing $x$, such that $Q$ stabilizes $U$ and acts continuously on the discrete ring $\Gamma(U,\Ow_X)$. We may also refer to a scheme endowed with an admissible right action of $Q$ as a \textit{$Q$-scheme}.

\subsection{\label{3.5}} Let $X$ be a scheme endowed with an admissible right action of a profinite group $Q$. Let $\Lambda$ be a noetherian ring. A \textit{$Q$-equivariant \'etale sheaf} of $\Lambda$-modules on $X$ is an \'etale sheaf of $\Lambda$-modules $\F$ on $X$ endowed with a collection $(\sigma(q))_{q \in Q}$ of isomorphisms $\F \rightarrow q_* \F$ such that:
\begin{enumerate}
\item we have $(q_{2*} \sigma(q_1)) \sigma(q_2) = \sigma(q_1q_2)$ for any $q_1,q_2$ in $Q$;
\item for any $Q$-equivariant \'etale morphism $U \rightarrow X$, with $U$ quasi-compact, the left action of $Q$ on the discrete set $\F(U)$ is continuous.
\end{enumerate}

We denote by $\Sh(X,Q,\Lambda)$ the corresponding category, by $D(X,Q,\Lambda)$ its the derived category, by $D_c^b(X,Q,\Lambda)$ the corresponding category of bounded complexes with constructible cohomology, and by $D_{ctf}^b(X,Q,\Lambda)$ the subcategory of complexes locally of finite Tor-amplitude.

\subsection{\label{3.22}} Let $\Lambda$ be a noetherian ring of finite global dimension. For any integer $n \geq 1$, we denote by $\Lambda_n$ the ring $\Lambda/n\Lambda$. Let $\widehat{\Lambda}$ be the projective limit of $(\Lambda_n)_{n \geq 1}$. Let $X$ be a scheme endowed with an admissible right action of a profinite group $Q$. A \textit{$Q$-equivariant $\widehat{\Lambda}$-system} on $X$ is an object $\F = (\F_n)_{n \geq 1}$ of the category
$$
\prod_{n \geq 1} D_{ctf}^b \left( X[\frac{1}{n}],Q,\Lambda_n \right),
$$
with isomorphisms $\theta_{n,m} : \F_{m} \otimes^{L}_{\Lambda_{m}} \Lambda_n \rightarrow (\F_n)_{|X[\frac{1}{m}]} $ for $n$ dividing $m$, such that $\theta_{l,n} \theta_{n,m} = \theta_{l,m}$. The $Q$-equivariant $\widehat{\Lambda}$-systems on $X$ form a triangulated category, which we abusively denote by $D_{c}^b(X,Q,\widehat{\Lambda})$.

\subsection{\label{3.24}} Let $X$ be a quasicompact quasiseparated scheme endowed with an admissible right action of a profinite group $Q$. We define a category $\widetilde{M}^{\mathrm{eff}}(X,Q,\Ow_E)$ as the smallest triangulated subcategory of $D_{c}^b(X,Q,\widehat{\Ow}_E)$ (cf. \ref{3.22}) which contains all $Q$-equivariant $\widehat{\Ow_E}$-systems of the form
$$
(Ra_! ( \F \otimes^L_{\Ow_E} \Ow_E/n\Ow_E) )_{n \geq 1},
$$
where $a : Z \rightarrow X$ is a separated morphism of finite presentation and $\F$ is in $D_{c}^b(Z,Q,\Ow_E)$ (cf. \ref{3.5}).

\subsection{\label{3.23}} Let $X$ be a scheme endowed with an admissible right action of a profinite group $Q$. The category $M^{\mathrm{eff}}(X,Q,\Ow_E)$ of \textit{effective $Q$-equivariant $\Ow_E$-systems} on $X$ is the full subcategory of $D_{c}^b(X,Q,\widehat{\Ow_E})$ (cf. \ref{3.22}) consisting of $Q$-equivariant $\widehat{\Ow_E}$-systems $\F$ on $X$ with the following property: for any point $x$ of $X$, there exists a $Q$-invariant affine open subscheme $U$ such that $\F_{|U} = (\F_{n|U})_n$ belongs to $\widetilde{M}(X,Q,\Ow_E)$ (cf. \ref{3.24}).

\begin{prop}\label{prop3.23.1} Let $X$ be a quasicompact quasiseparated scheme endowed with an admissible right action of a profinite group $Q$. Then the inclusion functor
$$
\widetilde{M}^{\mathrm{eff}}(X,Q,\Ow_E) \rightarrow M^{\mathrm{eff}}(X,Q,\Ow_E),
$$
is an equivalence of categories.
\end{prop}

We first remark that if $j : U \rightarrow X$ is the immersion of a $Q$-invariant quasicompact open subscheme of $X$, then the functor $j_! : \F \mapsto (j_! \F_n)_n$ from $D_{c}^b(U,Q,\widehat{\Ow}_E)$ to $D_{c}^b(X,Q,\widehat{\Ow}_E)$ sends $\widetilde{M}^{\mathrm{eff}}(U,Q,\Ow_E)$ to $\widetilde{M}^{\mathrm{eff}}(X,Q,\Ow_E)$. Indeed, if $a : Z \rightarrow U$ is a separated morphism of finite presentation and if $\F$ is an object of $D_{c}^b(Z,Q,\Ow_E)$, then we have
$$
j_! (Ra_! ( \F \otimes^L_{\Ow_E} \Ow_E/n\Ow_E) )_n \simeq (R(aj)_! ( \F \otimes^L_{\Ow_E} \Ow_E/n\Ow_E) )_n,
$$
and $aj : U \rightarrow X$ is a separated $Q$-equivariant morphism of finite presentation.

Let us now consider an object $\F$ of $M^{\mathrm{eff}}(X,Q,\Ow_E)$. Since $X$ is quasicompact, there exists a finite cover $(U_r)_{r \in R}$ by $Q$-invariant quasicompact open subschemes of $X$ such that $\F_{|U_r}$ belongs to $\widetilde{M}^{\mathrm{eff}}(U_r,Q,\Ow_E)$ for each $r$ in $R$. We prove by induction on the cardinality of $R$ that $\F$ belongs to $\widetilde{M}^{\mathrm{eff}}(X,Q,\Ow_E)$, the case where $R$ has at most one element being trivial. We denote by $j_U : U \rightarrow X$ the open immersion associated to an open subset $U$ of $X$. Let $r_0$ be an element of $R$, and let us consider $U=U_{r_0}$ and $V = \bigcup_{r \neq r_0} U_r$. Then we have a distinguished triangle
$$
(j_{U \cap V})_! j_{U \cap V}^* \F \rightarrow (j_U)_! j_U^* \F \oplus (j_V)_! j_V^* \F \rightarrow \F \xrightarrow[]{[1]}.
$$
The objects $(j_U)_! j_U^* \F, (j_V)_!j_V^*  \F$ and $(j_{U \cap V})_! j_{U \cap V}^*  \F$ of $D_c^b(X,Q,\widehat{\Ow}_E)$ belong to $\widetilde{M}^{\mathrm{eff}}(X,Q,\Ow_E)$ by induction and by the first part of the proof, hence the conclusion.

\begin{prop}\label{prop3.23.2} Let $Q$ be a profinite group, and let $f : X \rightarrow Y$ be a $Q$-equivariant morphism of schemes.
\begin{enumerate}
\item The functor $f^* : \F \mapsto (f^* \F_n)_n$ from $D_{c}^b(Y,Q,\widehat{\Ow}_E)$ to $D_{c}^b(X,Q,\widehat{\Ow}_E)$ sends $M^{\mathrm{eff}}(Y,Q,\Ow_E)$ to $M^{\mathrm{eff}}(X,Q,\Ow_E)$.
\item If $f$ is separated of finite presentation, then the functor $Rf_! : \F \mapsto (Rf_! \F_n)_n$ from $D_{c}^b(X,Q,\widehat{\Ow}_E)$ to $D_{c}^b(Y,Q,\widehat{\Ow}_E)$ sends $M^{\mathrm{eff}}(X,Q,\Ow_E)$ to $M^{\mathrm{eff}}(Y,Q,\Ow_E)$.
\end{enumerate}

\end{prop}

For the first part it is enough to prove that for any affine $Q$-invariant open subschemes $U$ and $V$ of $X$ and $Y$ respectively, with $f(U)$ contained in $V$, then the functor $f^*$ from $D_{c}^b(V,Q,\widehat{\Ow}_E)$ to $D_{c}^b(U,Q,\widehat{\Ow}_E)$ sends $\widetilde{M}(V,Q,\Ow_E)$ to $\widetilde{M}(U,Q,\Ow_E)$. This follows from the following remark: if $a : Z \rightarrow V$ is a separated morphism of finite presentation and if $\F$ is an object of $D_{c}^b(Z,Q,\Ow_E)$, then by the proper base change theorem we have
$$
f^* (Ra_! ( \F \otimes^L_{\Ow_E} \Ow_E/n\Ow_E) )_n \simeq (Ra'_! ( (f')^* \F \otimes^L_{\Ow_E} \Ow_E/n\Ow_E) )_n,
$$
where $a' : Z \times_V U \rightarrow U$ is the base change of $a$ along $f : U \rightarrow V$, hence is a separated $Q$-equivariant morphism of finite presentation, and $f' : Z \times_V U \rightarrow Z$ is the first projection.

We now prove $(2)$. Let $\F$ be an object of $M^{\mathrm{eff}}(X,Q,\Ow_E)$. Let $V$ be any affine $Q$-invariant open subscheme of $Y$. Then $(Rf_! \F)_{|V}$ is isomorphic to $R(f_{|U})_! \F_{|U}$, where $U = f^{-1}(V)$. Since $f$ is quasicompact and quasiseparated, we obtain that $U$ is quasicompact and quasiseparated. By Proposition \ref{prop3.23.1}, the object $\F_{|U}$ of $D_{c}^b(U,Q,\widehat{\Ow}_E)$ belongs to $\widetilde{M}^{\mathrm{eff}}(U,Q,\Ow_E)$. It remains to prove that the functor $R(f_{|U})_! $ from $D_{c}^b(U,Q,\widehat{\Ow}_E)$ to $D_{c}^b(V,Q,\widehat{\Ow}_E)$ sends $\widetilde{M}^{\mathrm{eff}}(U,Q,\Ow_E)$ to $M^{\mathrm{eff}}(V,Q,\Ow_E)$. This follows from the following remark: if $a : Z \rightarrow U$ is a separated morphism of finite presentation and if $\F$ is an object of $D_{c}^b(Z,Q,\Ow_E)$, then we have
$$
R(f_{|U})_!  (Ra_! ( \F \otimes^L_{\Ow_E} \Ow_E/n\Ow_E) )_n \simeq (R(f_{|U}a)_! (  \F \otimes^L_{\Ow_E} \Ow_E/n\Ow_E) )_n,
$$
and the composition $f_{|U}a : Z \rightarrow V$ is a separated $Q$-equivariant morphism of finite presentation.

\begin{cor}\label{cor3.23.20} Let $Q$ be a profinite group, and let $X$ be a $Q$-scheme. The Tate twist functor $\F \mapsto \F(-1)$ from $D_{c}^b(X,Q,\widehat{\Ow}_E)$ to $D_{c}^b(X,Q,\widehat{\Ow}_E)$ sends $M^{\mathrm{eff}}(X,Q,\Ow_E)$ to itself.
\end{cor}

Indeed, we have $\F(-1) \simeq Ra_! a^* \F[2]$, where $a : \mathbb{A}^1_X \rightarrow X$ is the relative affine line, and the conclusion then follows from Proposition \ref{prop3.23.2}.

\subsection{\label{3.25}} Let $X$ be a scheme endowed with an admissible right action of a profinite group $Q$. We have a triangulated functor
$$
D_{c}^b(X,Q,\Ow_E) \rightarrow M^{\mathrm{eff}}(X,Q,\Ow_E),
$$
which sends an object $\F$ of $D_{c}^b(X,Q,\Ow_E)$ to $(\F \otimes^L_{\Ow_E} \Ow_E/n\Ow_E)_n$. We simply denote by $\F \mapsto \F$ this functor.

\begin{prop}\label{prop3.23.5} Let $Q$ be a profinite group, let $(X_r)_{r \in R}$ be a $Q$-equivariant filtered diagram of quasicompact quasiseparated schemes with affine transition morphisms, and let $X$ be the limit of this diagram. Then the essential image of the canonical functor
$$
2-\mathrm{colim}_{r} M^{\mathrm{eff}}(X_r,Q,\Ow_E) \rightarrow M^{\mathrm{eff}}(X,Q,\Ow_E),
$$
generates $M^{\mathrm{eff}}(X,Q,\Ow_E)$ as a triangulated category.
\end{prop}

Let us consider an object $\F$ of $M^{\mathrm{eff}}(X,Q,\Ow_E)$ of the form $Ra_! \G$, where $a : Z \rightarrow X$ is a $Q$-equivariant separated morphism of finite presentation and $\G$ is an object of $D_{c}^b(Z,\Ow_E)$. Then there exists an index $r_1$ in $R$ such that $a$ is the pullback along $X \rightarrow X_{r_1}$ of a $Q$-equivariant separated morphism of finite presentation $a_1 : Z_1 \rightarrow X_{r_1}$. The canonical functor
$$
2-\mathrm{colim}_{r \rightarrow r_1} D_{c}^b(Z_1 \times_{X_{r_1}} X_r,\Ow_E) \rightarrow D_{c}^b(Z,Q,\Ow_E),
$$
is an equivalence of categories, hence there exists some index $r_2 \rightarrow r_1$ in $R$ such that $\G$ is the pullback along $Z \rightarrow Z_1 \times_{X_{r_1}} X_{r_2}$ of an object of $D_{c}^b(Z_1 \times_{X_{r_1}} X_{r_2},\Ow_E)$, and the conclusion then follows by Proposition \ref{prop3.23.1}.

\begin{prop}\label{prop3.23.6} Let $Q$ be a profinite group, let $X$ be a quasicompact quasiseparated $Q$-scheme. Then the triangulated category $M^{\mathrm{eff}}(X,Q,\Ow_E)$ is generated by objects of the form $Ra_! \G$, where $a : Z \rightarrow X$ is a separated morphism of finite presentation and $\G$ is a $Q$-equivariant lisse \'etale sheaf of $\Ow_E$-modules on $Z$.
\end{prop}

By Proposition \ref{prop3.23.5}, we can assume (and we do) that $X$ is of finite type over $\mathbb{Z}$. Let us consider an object $\F$ of $M^{\mathrm{eff}}(X,Q,\Ow_E)$ of the form $Ra_! \G$, where $a : Z \rightarrow X$ is a $Q$-equivariant separated morphism of finite presentation and $\G$ is an object of $D_{c}^b(Z,\Ow_E)$. We argue by noetherian induction on $Z$. Let $j : U \rightarrow Z$ be a non empty open subset of $Z$ such that $j^*\G$ has lisse cohomology sheaves, and let $i : T \rightarrow Z$ be a closed immersion with support $Z \setminus U$. Then we have a distinguished triangle
$$
R(aj)_! j^* \G \rightarrow Ra_! \G \rightarrow R(ai)_!i^* \G \xrightarrow[]{[1]}.
$$
The conclusion sought holds for $R(ai)_!i^* \G$ by noetherian induction, and it clearly holds for $R(aj)_! j^* \G$ as well, hence the result.

\begin{prop}\label{prop3.23.8} Let $Q$ be a profinite group and let $f : X \rightarrow Y$ be a $Q$-equivariant universal homeomorphism of $Q$-schemes. Then the equivalence of categories $f^*$ from $D_{c}^b(Y,Q,\widehat{\Ow}_E)$ to $D_{c}^b(X,Q,\widehat{\Ow}_E)$ induces an equivalence of categories from $M^{\mathrm{eff}}(Y,Q,\Ow_E)$ to $M^{\mathrm{eff}}(X,Q,\Ow_E)$.
\end{prop}

The functor $f^*$ sends $M^{\mathrm{eff}}(Y,Q,\Ow_E)$ to $M^{\mathrm{eff}}(X,Q,\Ow_E)$ by Proposition \ref{prop3.23.2}(1). It remains to prove that the quasi-inverse $f_*$ sends $M^{\mathrm{eff}}(X,Q,\Ow_E)$ to $M^{\mathrm{eff}}(Y,Q,\Ow_E)$. We can assume (and we do) that $Y = \Spec(A)$ is affine. Then the seminormalization $Y' = \Spec(A')$ of $Y^{\mathrm{red}}$ is a $Q$-scheme, and is a filtered limit of finitely presented $Q$-equivariant universal homeomorphisms to $Y$. Let us write $X' = (Y' \times_Y X)^{\mathrm{red}}$ as $\Spec(B')$. For any $Q$-equivariant factorization $A' \rightarrow C \rightarrow B'$ with $C$ finitely presented over $A'$, there exists by \cite[Cor. 23]{Kelly} a finitely generated ideal $I$ of $C$ contained in the kernel of the homomorphism $C \rightarrow B'$ such that $A \rightarrow C/I$ is a universal homeomorphism. By replacing $I$ with $\sum_{q \in Q} qI$ if necessary, the ideal $I$ can be taken to be $Q$-invariant. Thus $X'$ is a filtered limit of finitely presented $Q$-equivariant universal homeomorphisms to $Y'$. 

Thus there exists a filtered system $(Y_r)_{r \in R}$ of finitely presented $Q$-equivariant universal homeomorphisms to $Y$, and a $Q$-equivariant universal homeomorphism $\lim_r Y_r \rightarrow X$ of $X$-schemes. By Proposition \ref{prop3.23.5}, it is enough to prove that $f_*$ sends $M^{\mathrm{eff}}(X,Q,\Ow_E)$ to $M^{\mathrm{eff}}(Y,Q,\Ow_E)$ when $f$ is a finitely presented $Q$-equivariant universal homeomorphisms. Such a morphism $f$ is a finite morphism of finite presentation, hence $f_* = Rf_!$ sends $M^{\mathrm{eff}}(X,Q,\Ow_E)$ to $M^{\mathrm{eff}}(Y,Q,\Ow_E)$ by Proposition \ref{prop3.23.2}(2); this concludes the proof of Proposition $\ref{prop3.23.8}$.

\subsection{\label{3.26}} Given a system $(\F_n)_n$ in $M^{\mathrm{eff}}(X,Q,\Ow_E)$, with $X$ quasicompact quasiseparated, each $\F_n$ has Tor-amplitude in a finite interval $[a,b]$, and is given by a perfect complex on $U[\frac{1}{n}]$ for some dense open subset $U$ of $X$. We now deduce from the uniform constructibility results in \cite{Ill10} that the open subset $U$ and the Tor-amplitude $[a,b]$ can be taken uniformly in $n$:

\begin{prop}\label{prop3.23.3}  Let $X$ be a quasicompact quasiseparated scheme endowed with an admissible right action of a profinite group $Q$. For any object $\F$ of $M^{\mathrm{eff}}(X,Q,\Ow_E)$, the following assertions hold.
\begin{enumerate}
\item There exists relative integers $a \leq b$ such that for each integer $n \geq 1$, the object $\F_n$ of $D_{ctf}^b(X[\frac{1}{n}],Q,\Ow_E/n \Ow_E)$ has Tor-amplitude in $[a,b]$.
\item There exists a dense $Q$-invariant open subset $U$ of $X$ such that for each integer $n \geq 1$, the object $\F_{n|U}$ of $D_{ctf}^b(U[\frac{1}{n}],Q,\Ow_E/n \Ow_E)$ is perfect, i.e. its cohomology sheaves $\mathcal{H}^q(\F_{n|U})$ are lisse \'etale sheaves of $\Ow_E/n \Ow_E$-modules on $U[\frac{1}{n}]$.
\end{enumerate}
\end{prop}

We can assume (and we do) that $Q = 1$ is the trivial group, and, by Proposition \ref{prop3.23.5}, that $X$ is of finite type over $\mathbb{Z}$. By Proposition \ref{prop3.23.6}, we can further assume (and we do) that $\F$ is of the form $Ra_! \G$, where $a : Z \rightarrow X$ is a separated morphism of finite presentation and $\G$ is a lisse \'etale sheaf of $\Ow_E$-modules on $Z$. We can also assume (and we do) that $Z$ is normal and connected.

We proceed by noetherian induction on $Z$. Let $\overline{z}$ be a geometric point of $Z$. The \'etale fundamental group $\pi_1(Z,\overline{z})$ acts continuously on $\G_{\overline{z}}$, hence its action factors through a finite quotient $I$, corresponding to a finite \'etale Galois cover $h :\widetilde{Z} \rightarrow Z$. We then have
$$
\G \simeq (h_* h^* \G)^I \simeq (\G_{\overline{z}} \otimes_{\Ow_E}h_* \Ow_E)^I,
$$
as \'etale sheaves of $\Ow_E$-modules on $Z$, hence
$$
Ra_! \G_n \simeq R\Gamma^I( \G_{\overline{z}} \otimes_{\Ow_E}^L R(ah)_! (\Ow_E/n \Ow_E) ).
$$
Let $d$ be an integer such that all the fibers of $a$ have dimension at most $d$, and let $\overline{Z}$ be the schematic image of $a$. By \cite[Th. 2.1]{Ill10}, there exists a dense open subscheme $U$ of $\overline{Z}$ such that for any $n \geq 1$, the object $R(ah)_! (\Ow_E/n \Ow_E)_{|U}$ of $D_{ctf}^b(U[\frac{1}{n}],\Ow_E/n \Ow_E)$ is perfect, with cohomological amplitude in $[0,2d]$. Since $R(ah)_! (\Ow_E/n \Ow_E)_{|U} = Ra_!(h_* (\Ow_E/n \Ow_E))_{|U}$ is also in $D_{ctf}^b(U[\frac{1}{n}],\Ow_E/n \Ow_E[I])$, we deduce that it is perfect as an object of the latter category. Thus $(Ra_! \G_n)_{|U}$ is perfect, with cohomological amplitude in $[0,2d]$. The conclusion then follows by applying the induction hypothesis to the complement of $a^{-1}(U)$ in $Z$.

\subsection{\label{3.16}} Let $\varphi : Q' \rightarrow Q$ be a continuous homomorphism of profinite groups and let $X$ be a scheme endowed with an admissible right action of $Q$. Then we have a functor $\varphi^*$ from $D_{c}^b(X,Q,\widehat{\Ow}_E)$ to $D_{c}^b(X,Q',\widehat{\Ow}_E)$ obtained by letting $Q'$ act on $Q$-equivariant objects through $\varphi$.

\begin{prop} Let $\varphi : Q' \rightarrow Q$ be a continuous homomorphism of profinite groups and let $X$ be a scheme endowed with an admissible right action of $Q$. Then the functor $\varphi^*$ from $D_{c}^b(X,Q,\widehat{\Ow}_E)$ to $D_{c}^b(X,Q',\widehat{\Ow}_E)$ sends $M^{\mathrm{eff}}(X,Q,\Ow_E)$ to $M^{\mathrm{eff}}(X,Q',\Ow_E)$.
\end{prop}

Indeed, we can assume that $X$ is affine, in which case it is enough to check that for any separated morphism of finite presentation $a : Z \rightarrow X$ and any object $\F$ of $D_{c}^b(Z,Q,\Ow_E)$, the object $\varphi^* Ra_! \F$ of $D_{c}^b(X,Q',\widehat{\Ow}_E)$ belongs to $M^{\mathrm{eff}}(X,Q',\Ow_E)$. This is immediate since $\varphi^* Ra_! \F$ is simply $Ra_! \F$, with $a : Z \rightarrow X$ a $Q'$-equivariant separated morphism of finite presentation.

\begin{prop}\label{prop3.16.9} Let $X$ be a quasicompact quasiseparated scheme endowed with an admissible right action of a profinite group $Q$, and let $I$ be a closed subgroup of $Q$ acting trivially on $X$. Then the essential image of the canonical functor
$$
2-\mathrm{colim}_{I'} M^{\mathrm{eff}}(X,Q/I',\Ow_E) \rightarrow M^{\mathrm{eff}}(X,Q,\Ow_E),
$$
where $I'$ runs over open subgroups of $I$ which are normal in $Q$, generates $M^{\mathrm{eff}}(X,Q,\Ow_E)$ as a triangulated category.
\end{prop}

Let us consider an object $\F$ of $M^{\mathrm{eff}}(X,Q,\Ow_E)$ of the form $Ra_! \G$, where $a : Z \rightarrow X$ is a separated morphism of finite presentation and $\G$ is an object of $D_{c}^b(Z,\Ow_E)$. Then there exists an open subgroup $I'$ of $I$, normal in $Q$, which acts trivially on $(Z,\F)$, hence the conclusion.

\subsection{\label{3.6}} Let $\Lambda$ be a noetherian ring and let $X$ be a scheme endowed with an admissible right action of a profinite group $Q$. For any open subgroup $H$ of $Q$ and any $H$-equivariant \'etale $\Lambda$-sheaf $\F$ on $X$, we denote by $\Ind_H^Q \F$ the $Q$-equivariant \'etale $\Lambda$-sheaf on $X$ whose local sections are given by the collections $(s_q)_{q \in Q}$, where $s_q$ is a local section of $q^* \F$ with $h \cdot s_{hq} = s_q$ for any element $h$ of $H$. The functor $\Ind_H^Q$ extends to a triangulated functor
\begin{align*}
\Ind_H^Q : D_{ctf}^b(X,H,\Lambda)  \rightarrow D_{ctf}^b(X,Q,\Lambda).
\end{align*}
In particular, we obtain a triangulated functor
$$
\Ind_H^Q : D_{c}^b(X,H,\widehat{\Ow}_E) \rightarrow D_{c}^b(X,Q,\widehat{\Ow}_E).
$$

\begin{prop}\label{prop3.6.1.1} Let $X$ be a scheme endowed with an admissible right action of a profinite group $Q$, and let $H$ be an open subgroup of $Q$. Then the functor $\Ind_H^Q$ from $D_{c}^b(X,H,\widehat{\Ow}_E)$ to $D_{c}^b(X,Q,\widehat{\Ow}_E)$ sends $M^{\mathrm{eff}}(X,H,\Ow_E)$ to $M^{\mathrm{eff}}(X,Q,\Ow_E)$.
\end{prop}

We can assume (and we do) that $X$ is affine. Let us consider an object $\F$ of $M^{\mathrm{eff}}(X,Q,\Ow_E)$ of the form $Ra_! \G$, where $a : Z \rightarrow X$ is an $H$-equivariant separated morphism of finite presentation and $\G$ is an object of $D_{c}^b(Z,H,\Ow_E)$. Let us consider the morphism $a' : Z' \rightarrow X$ defined by
$$
a' : Z' = \colim_{q \in Q} Z \times_{X,q} X \rightarrow X,
$$
with transition isomorphisms $Z \times_{X,q} X \rightarrow Z \times_{X,q h} X$ given by $h \times \id_X$ for $h$ in $H$. The $X$-scheme $Z'$ is separated of finite presentation since $H$ has finite index in $Q$, and is endowed with an admissible action of $Q$ by setting the action of an element $t$ of $Q$ to be the automorphism induced by
$$
Z \times_{X,q} X \xrightarrow[]{\id_Z \times t } Z \times_{X,t^{-1}q} X.
$$
If $\G'$ is the pullback of $\G$ to $Z'$ by the first projection, then we have
$$
Ra'_! \G' \simeq \Ind_{H}^Q Ra_! \G ,
$$
hence the conclusion.

\subsection{\label{3.19}} Let $X$ be a scheme endowed with an admissible right action of a profinite group $Q$. We define a biadditive pairing
\begin{align*}
\otimes_{\Ow_E}^L :  D_{c}^b(X,Q,\widehat{\Ow}_E) \times  D_{c}^b(X,Q,\widehat{\Ow}_E) \rightarrow  D_{c}^b(X,Q,\widehat{\Ow}_E)
\end{align*}
by sending $(\F,\G)$ to $(\F_n \otimes^L_{\Ow_E/n\Ow_E} \G_n)_n$.

\begin{prop}\label{prop3.19.10} Let $X$ be a scheme endowed with an admissible right action of a profinite group $Q$. Then the functor $\otimes_{\Ow_E}^L$ from $D_{c}^b(X,Q,\widehat{\Ow}_E) \times  D_{c}^b(X,Q,\widehat{\Ow}_E)$ to $D_{c}^b(X,Q,\widehat{\Ow}_E)$ sends $M^{\mathrm{eff}}(X,Q,\Ow_E) \times M^{\mathrm{eff}}(X,Q,\Ow_E)$ to $M^{\mathrm{eff}}(X,Q,\Ow_E)$.
\end{prop}

We can assume (and we do) that $X$ is affine, and we use Proposition \ref{prop3.23.1}. Let us consider $Q$-invariant separated morphisms of finite presentation $a : Z \rightarrow X$ and $a' : Z' \rightarrow X$. Let $\F$ be in $D_c^b(Z,Q,\Ow_E)$ and let $\F'$ be in $D_{c}^b(Z',Q,\Ow_E)$. By the K\"unneth formula we have
$$
Ra_! \F \otimes_{\Ow_E}^L Ra'_! \G \simeq Rb_! (\mathrm{pr}_1^*\F \otimes_{\Ow_E}^L \mathrm{pr}_2^*\G),
$$
where $b : Z \times_X Z' \rightarrow X$ is the fiber product and $\mathrm{pr}_1,\mathrm{pr}_2$ are the first and second projections on $ Z \times_X Z'$. The morphism $b$ is separated of finite presentation, and $\mathrm{pr}_1^*\F \otimes_{\Ow_E}^L \mathrm{pr}_2^*\G$ is in $D_{c}^b(Z \times_X Z',Q,\Ow_E)$.

\subsection{\label{3.27}} Let $X$ be a scheme endowed with an admissible right action of a profinite group $Q$. We define a category $M^{\mathrm{eff}}(X,Q,E)$ of \textit{effective $Q$-equivariant $E$-systems}, with the same objects as $M^{\mathrm{eff}}(X,Q,\Ow_E)$ and with morphisms defined by
$$
\Hom_{M^{\mathrm{eff}}(X,Q,E)}(\F,\G) = \mathrm{lim}_{U} E \otimes_{\Ow_E} \Hom_{M^{\mathrm{eff}}(U,Q,\Ow_E)}(\F_{|U},\G_{|U})  ,
$$
for $\F,\G$ in $M^{\mathrm{eff}}(X,Q,\Ow_E)$, where $U$ runs over quasicompact quasiseparated open subschemes of $X$. We similarly define a category $D_{c}^b(X,Q,\widehat{E})$ by replacing $M^{\mathrm{eff}}(X,Q,\Ow_E)$
with $D_{c}^b(X,Q,\widehat{\Ow}_E)$. If we replace $M^{\mathrm{eff}}(X,Q,\Ow_E)$ by $D_{c}^b(X,Q,\Ow_E)$ in this definition, we simply obtain $D_c^b(X,Q,E)$. Correspondingly, we have a natural functor
$$
D_c^b(X,Q,E) \rightarrow M^{\mathrm{eff}}(X,Q,E),
$$
which we simply denote by $\F \mapsto \F$.

\begin{prop}\label{prop3.26.10} Let $X$ be a scheme endowed with an admissible right action of a profinite group $Q$, and let $I$ be a finite normal subgroup of $Q$ acting trivially on $X$. Then the functor $R \Gamma^I$ from $D_{c}^b(X,Q,\widehat{E})$ to $D_{c}^b(X,Q/I,\widehat{E})$ sends $M^{\mathrm{eff}}(X,Q,E)$ to $M^{\mathrm{eff}}(X,Q/I,E)$.
\end{prop}

We can assume (and we do) that $X$ is affine. Let us consider an object $\F$ of $M^{\mathrm{eff}}(X,Q,\Ow_E)$ of the form $Ra_! \G$, where $a : Z \rightarrow X$ is a $Q$-equivariant separated morphism of finite presentation and $\G$ is an object of $D_{c}^b(Z,Q,E)$. Let us consider the factorisation
$$
Z \xrightarrow[]{c} Z/ I \xrightarrow[]{b} X,
$$
of $a$. Then $Z/ I$ is endowed with an admissible action of $Q/ I$, the morphism $b$ is $Q/I$-equivariant, and since $c$ is finite, we have
$$
R \Gamma^I( Ra_! \G) \simeq Rb_! R \Gamma^I(c_*\G) ,
$$
in $D_{c}^b(X,Q/I,\widehat{E})$, and $R \Gamma^I(c_*\G)$ is the image of an object of $D_{c}^b(Z,Q,E)$.

\begin{prop}\label{prop3.19.1} Let $X$ be a quasicompact quasiseparated scheme endowed with an admissible right action of a profinite group $Q$. Then the triangulated category $M^{\mathrm{eff}}(X,Q,E)$ is generated by $Q$-equivariant $E$-systems of the form
$$
R\Gamma^I(V \otimes_E Ra_! E),
$$
where $a : Z \rightarrow X$ is a $Q'$-equivariant separated morphism of finite presentation, for some extension $Q'$ of $Q$ by a finite group $I$, and where $V$ is a finite dimensional continuous $E$-representation of $Q'$. 
\end{prop}

By Proposition \ref{prop3.23.5}, we can assume (and we do) that $X$ is of finite type over $\mathbb{Z}$. We then use Proposition \ref{prop3.23.6}; let us consider an object $\F$ of $M^{\mathrm{eff}}(X,Q,E)$ of the form $Ra_! \G$, where $a : Z \rightarrow X$ is a $Q$-equivariant separated morphism of finite type and $\G$ is a lisse \'etale sheaf of $E$-modules on $Z$. We can further assume (and we do) that $Z$ is normal and connected.

Let $\overline{z}$ be a geometric point, and let $\mathrm{ev}_{\overline{z}}$ be the functor from the category of finite \'etale $Z$-schemes to the category of sets, sending such a $Z$-scheme $Y$ to the finite set $Y(\overline{z})$. Let $\tilde{Q}$ be the profinite group of pairs $(c,q)$, where $q$ is an element of $Q$ and $c : \mathrm{ev}_{\overline{z}} \rightarrow \mathrm{ev}_{\overline{z} q}$ is an isomorphism of functors, with composition law $(c,q)(c',q') = ((c(q')^*) \circ c', qq')$. The profinite group $\tilde{Q}$ is an extension of $Q$ by $\pi_1(Z,\overline{z})$. Then $\G_{\overline{z}}$ is a continuous finite dimensional $E$-representation of $\tilde{Q}$. Let $H$ be an open subgroup of $\pi_1(Z,\overline{z})$, normal in $\tilde{Q}$, such that $H$ acts trivially on $\G_{\overline{z}}$. Then $\G_{\overline{z}}$ is a representation of $Q' = \tilde{Q}/H$, which is an extension of $Q$ by the finite group $I = \pi_1(Z,\overline{z})/H$. Let $b : Z' \rightarrow Z$ be a finite \'etale $Q'$-equivariant morphism, which is Galois of group $I$. We then have
$$
Ra_! \G \simeq R\Gamma^I(\F_{\overline{z}} \otimes_E R(ab)_! E),
$$
hence the conclusion.

\begin{cor}\label{cor3.19.2} Let $X$ be a quasicompact quasiseparated scheme endowed with an admissible right action of a profinite group $Q$. Then the triangulated category $M^{\mathrm{eff}}(X,Q,E)$ is generated by $Q$-equivariant $E$-systems of the form
$$
R\Gamma^I(V \otimes_E Ra_* E),
$$
where $a : Z \rightarrow X$ is a $Q'$-equivariant proper morphism of finite presentation, for some extension $Q'$ of $Q$ by a finite group $I$, and where $V$ is a finite dimensional continuous $E$-representation of $Q'$. 
\end{cor}

\section{Vanishing and nearby cycles \label{nearby}}

We fix a number field $E$, with ring of integers $\Ow_E$.

\subsection{\label{3.7}} Let $S$ be the spectrum of a henselian valuation ring with separably closed fraction field, endowed with an admissible right action from a profinite group $Q$. We say that $S$ is a $Q$-trait if the quotient $S/Q$ is the spectrum of a (necessarily henselian) discrete valuation ring. We say that the $Q$-trait $S$ is $Q$-excellent if $S/Q$ is excellent.

\subsection{\label{3.8}} Let $S$ be a $Q$-trait, cf. \ref{3.7}, such that the kernel of the homomorphism $Q \rightarrow \Aut(S)^{\mathrm{op}}$ is finite. Let $j : \eta \rightarrow S$ and be the generic point of $S$ and let $i : s \rightarrow S$. Let $S_0$ be the quotient trait $S/Q$, with generic point $\eta_0$ and special point $s_0$. Then $\eta$ is a separable closure of $\eta_0$, the $Q$-trait $S$ is the normalization of $S_0$ in $\eta$, and we have an exact sequence
$$
1 \rightarrow \Gamma \rightarrow Q \rightarrow \Gal(\eta/\eta_0) \rightarrow 1,
$$
for some finite group $\Gamma$. There exists a factorisation $\eta \rightarrow \eta' \rightarrow \eta_0$, with $\eta' \rightarrow \eta_0$ a Galois extension, such that the exact sequence above is the pullback along the surjection from $\Gal(\eta/\eta_0)$ to $\Gal(\eta'/\eta_0)$ of an exact sequence of finite groups
$$
1 \rightarrow \Gamma \rightarrow Q' \rightarrow \Gal(\eta'/\eta_0) \rightarrow 1.
$$
We thus have $Q = Q' \times_{\Gal(\eta'/\eta_0)} \Gal(\eta/\eta_0)$, hence an exact sequence
$$
1 \rightarrow \Gal(\eta/\eta') \rightarrow Q \rightarrow  Q' \rightarrow 1.
$$

\begin{lem}\label{lem3.8.1} Let $\Lambda$ be a ring. For any $Q'$-equivariant morphism $X \rightarrow \eta'$, the restriction functor
$$
\Sh(X,Q',\Lambda) \rightarrow \Sh(X_{\eta},Q,\Lambda), 
$$
is an equivalence of categories.
\end{lem}

Indeed, any object $\F$ of $\Sh(X_{\eta},Q,\Lambda)$ is equipped with a continuous action of the subgroup $\Gal(\eta/\eta')$ of $Q$ compatible with the Galois action on $X_{\eta} = X \times_{\eta'} \eta$, hence is canonically the pullback to $X_{\eta}$ of an \'etale sheaf of $\Lambda$-modules $\G$ on $X$. The action of $Q$ on $\F$ then translates to an action of $Q'$ on $\G$.

\subsection{\label{3.12}} Let $S$ be a $Q$-trait, cf. \ref{3.7}, and let $S_0$ be as in \ref{3.8}. We define the completion $\hat{S}$ as $S \times_{S_0} \hat{S_0}$, where $\hat{S_0}$ is the completion of the trait $S_0$. Then $\hat{S}$ is naturally a $Q$-trait, and is always $Q$-excellent. One should note that $\hat{S}$ typically differs from the completion of the local scheme $S$ at its closed point.

\subsection{\label{3.11}} Let $S$ be a $Q$-trait, cf. \ref{3.7}. Let $\Lambda$ be a finite ring of order invertible on $S$. Let us consider a $Q$-equivariant morphism $X \rightarrow S$ which is locally of finite presentation. We have a canonical $Q$-equivariant diagram
$$
X_{\eta} \xrightarrow[]{j_X} X \xleftarrow[]{i_X} X_s.
$$
Then we define a \textit{nearby cycle functor} as follows:
\begin{align*}
R\Psi_{X/S} : D^+(X_{\eta}, Q, \Lambda) &\rightarrow D^+(X_{s}, Q, \Lambda) \\
\F &\mapsto i_X^* Rj_{X*} \F.
\end{align*}

\begin{prop}\label{prop3.11.1} Let $S,Q,X, \Lambda$ be as in \ref{3.11}. Then the functor
$$
R\Psi_{X/S} : D^+(X_{\eta}, Q, \Lambda) \rightarrow D^+(X_{s}, Q, \Lambda),
$$
sends $D_{ctf}^b(X_{\eta}, Q, \Lambda) $ to $D_{ctf}^b(X_{s}, Q, \Lambda)$. Moreover, if $\hat{\eta}$ is the generic point of $\hat{S}$ (cf. \ref{3.12}), then $R\Psi_{X/S}$ coincides with the composition
$$
D^+(X_{\eta}, Q, \Lambda) \rightarrow D^+(X_{\hat{\eta}}, Q, \Lambda) \xrightarrow[]{R\Psi_{X_{\hat{S}}/\hat{S}}} D^+(X_{s}, Q, \Lambda)
$$
\end{prop}

The problem is local on $X$, hence we can assume that $X$ is of finite presentation over $S$. Since $\Lambda$ is finite, we can further assume (and we do) that the kernel of the homomorphism $Q \rightarrow \Aut(S)^{\mathrm{op}}$ is finite. Let $\eta \rightarrow \eta' \rightarrow \eta_0$ be as in \ref{3.8}. We keep the notation from \ref{3.8} and we consider the normalization $S'$ of $S_0$ in $\eta'$. By appropriately enlarging $\eta$, we can assume (and we do) that $X$ is of the form $X' \times_{S'} S$ for some finitely presented $Q'$-equivariant $S'$-scheme $X'$. Let $s'$ be the closed point of $S'$ and $\overline{s}'$ be its separable closure in $s$, so that $s$ is a purely inseparable extension of $\overline{s}'$. Then the composition (cf. \ref{lem3.8.1})
$$
D^+(X_{\eta'}, Q', \Lambda) \simeq D^+(X_{\eta}, Q, \Lambda) \xrightarrow[]{R\Psi_{X/S}} D^+(X_{s}, Q, \Lambda) \simeq D^+(X'_{\overline{s}'}, Q, \Lambda) \rightarrow D^+(X'_{\overline{s}'}, \Gal(\eta/\eta'), \Lambda)
$$
coincides with the usual nearby cycle functor for $X' \rightarrow S'$, which preserves constructibility by \cite[Th. Finitude, 3.2]{SGA412}, hence the first assertion of the Proposition. We then remark that $\hat{S}$ coincides with $S \times_{S'} \hat{S'}$, and that the composition
$$
D^+(X_{\eta'}, Q', \Lambda) \rightarrow D^+(X_{\hat{\eta'}}, Q', \Lambda) \xrightarrow[]{R \Psi_{X'_{\hat{S'}}/\hat{S'}}}  D^+(X'_{\overline{s}'}, Q, \Lambda),
$$
where $\hat{\eta'}$ is the generic point of $\hat{S'}$, coincides with the nearby cycle functor for $X' \rightarrow S'$ by \cite[Th. Finitude, 3.7]{SGA412}, hence the last assertion.

\begin{cor} Let $S,Q,X$ be as in \ref{3.11}. Then $R \Psi_{X/S}$ induces a functor from $D_{c}^b(X_{\eta},Q,\widehat{\Ow}_E)$ to $D_{c}^b(X_{s},Q,\widehat{\Ow}_E)$.
\end{cor}

\subsection{\label{3.3}} We are now ready to state the main result of this section, whose proof follows, as \cite[4]{Zh09} and \cite[5.2]{Ka19}, the method initiated by Vidal in \cite{Vi}.

\begin{teo}\label{teo3.3.1} Let $Q$ be a profinite group, let $S$ be a $Q$-trait, cf. \ref{3.7}, with generic point $\eta$ and closed point $s$. Let $X$ be a finitely presented $Q$-equivariant $S$-scheme. Then the functor $R \Psi_{X/S}$ from $D_{c}^b(X_{\eta},Q,\widehat{E})$ to $D_{c}^b(X_{s},Q,\widehat{E})$ sends $M^{\mathrm{eff}}(X_{\eta},Q,E)$ to $M^{\mathrm{eff}}(X_{s},Q,E)$.

\end{teo}

We follow closely the proof of \cite[5.2]{Ka19}. By Proposition \ref{prop3.19.1}, it is enough to prove that objects of the form
$$
R\Gamma^I(V \otimes_E Ra_* E),
$$
are sent by $R\Psi_{X/S}$ to $M^{\mathrm{eff}}(X_{s},Q,E)$, where $a : Z \rightarrow X$ is a proper $Q'$-equivariant separated morphism of finite presentation, for some extension $Q'$ of $Q$ by a finite group $I$, and where $V$ is a finite dimensional continuous $E$-representation of $Q'$. 

Let $\overline{a} : \overline{Z} \rightarrow X$ be a proper finitely presented $Q$-equivariant morphism whose pullback along $X_{\eta} \rightarrow X$ coincides with $a$. Let $\overline{a}_s : \overline{Z}_s \rightarrow X_s$ be the corresponding morphism on the special fibers. Then the proper base change theorem yields
$$
R \Psi_{X/S}  R\Gamma^I(V \otimes_E Ra_* E) \simeq R\Gamma^I(V \otimes_E R(\overline{a}_s)_* R \Psi_{\overline{Z}/S} E),
$$
hence Theorem \ref{teo3.3.1} follows from the following lemma, applied to the finitely presented $Q'$-equivariant $S$-scheme $\overline{Z}$.

\begin{lem}\label{lem3.3.2} Let $Q$ be a profinite group. Let $S$ be a $Q$-trait and let $X$ be a finitely presented $Q$-equivariant $S$-scheme. Then the object $R \Psi_{X/S} E$ of $D_{c}^b(X_{s},Q,\widehat{E})$ belongs to $M^{\mathrm{eff}}(X_{s},Q,E)$.
\end{lem}

We prove Lemma \ref{lem3.3.2} by induction on the dimension $d$ of $X_{\eta}$, the case $d < 0$ being trivial. We can assume (and we do) that $Q$ acts faithfully on the pair $(X,S)$, so that the kernel of the homomorphism $Q \rightarrow \Aut(S)^{\mathrm{op}}$ is finite. By Proposition \ref{prop3.11.1}, we can assume (and we henceforth do) that $S$ is $Q$-excellent. Let $\eta \rightarrow \eta' \rightarrow \eta_0$ be as in \ref{3.8}. Let $S'$ be the normalization of $S_0$ in $\eta'$. We can assume (and we do) that $X$ is of the form $X' \times_{S'} S$ for some finitely presented $Q'$-equivariant $S'$-scheme $X'$. Let $s'$ be the closed point of $S'$. We need to prove that the object $R\Psi_{X'/S'} E$ of $D_{c}^b(X_{s},Q,\widehat{E})$ belongs to $M^{\mathrm{eff}}(X_{s},Q,E)$. We proceed in several steps.

\begin{enumerate}
\item Let $\tilde{\eta} \rightarrow \eta$ be a finite field extension and let $\tilde{S}$ be the normalization of $S$ in $\tilde{\eta}$. Since $\eta$ is separably closed, this extension is purely inseparable, and there is a unique admissible action of $Q$ on $\tilde{\eta}$ such that $\tilde{\eta} \rightarrow \eta$ is $Q$-equivariant. By enlarging $\eta'$ if necessary, we can assume that $\tilde{\eta}$ is given by $\tilde{\eta}' \times_{\eta' } \eta$, for some finite purely inseparable extension $\tilde{\eta}' \rightarrow \eta'$. Let $\tilde{S}'$ be the normalization of $S'$ in $\tilde{\eta}'$. Then by \cite[Th. Finitude, 3.7]{SGA412}, we have
$$
(R\Psi_{X'/S'} \Ow_E/n\Ow_E)_{|X_s[\frac{1}{n}]} \simeq (R\Psi_{X'_{\tilde{S}'}/\tilde{S}'}  \Ow_E/n\Ow_E)_{|X_s[\frac{1}{n}]}.
$$
Consequently, the conclusion of Lemma \ref{lem3.3.2} holds for the $Q$-equivariant morphism $X \rightarrow S$ if and only if it holds for the $Q$-equivariant morphism $X_{\tilde{S}} \rightarrow\tilde{S}$.

\item By $(1)$ and by enlarging $\eta'$ if necessary, we can assume (and we do) that the irreducible components of $X'_{\eta'}$ are geometrically integral over $\eta'$. Let $a : X'' \rightarrow X'$ be the normalization of $X'$, which is an isomorphism above some open dense subscheme $j : U \rightarrow X'$. Let $Z'$ be the reduced complement of $U$ in $X'$ and $Z''$ its inverse image by $a$.
Since the generic fibers of $Z''$ and $Z'$ have dimension less than $d$, the induction hypothesis ensures that $R\Psi_{Z''/S'}E$ and $R\Psi_{Z'/S'} E$ belong to $M^{\mathrm{eff}}(X_{s},Q,E)$. Thus $R\Psi_{X'/S'} E$ is in $M^{\mathrm{eff}}(X_{s},Q,E)$, if and only if so is $R\Psi_{X'/S'} j_!E$, hence if and only if $R\Psi_{X''/S'} E$ is in $M^{\mathrm{eff}}(X_{s},Q,E)$.

\item By $(2)$, we can assume (and we do) that $X'$ is normal and that we have a finite decomposition 
$$
X' = \bigsqcup_{j \in J} Y^{(j)}
$$ 
into disjoint open subsets, where $Y^{(j)}$ is irreducible with $Y^{(j)}_{\eta'}$ geometrically integral over $\eta'$. The profinite group $Q$ acts continuously on the right on $J$, hence the stabilizer $Q_j$ of an element $j$ of $J$ is an open subgroup of $Q$. We have
$$
R\Psi_{X/S} E = \sum_{j \in J/Q} \Ind_{Q_j}^Q \tau^{(j)}_* R \Psi_{Y_S^{(j)}/S} E,
$$
in $D_{c}^b(X_{s},Q,\widehat{E})$, where $\tau^{(j)} : Y^{(j)}_s \rightarrow X_s$ is the natural closed immersion. By Proposition \ref{prop3.6.1.1}, if the conclusion of Lemma \ref{lem3.3.2} holds for the $Q_j$-equivariant morphism $Y_S^{(j)} \rightarrow S$ for each $j$ in $J$, then it holds for $X \rightarrow S$ as well.

\item By $(3)$, we can assume (and we do) that $X'$ is normal and that $X'_{\eta'}$ is geometrically integral. By \ref{3.16}, we can assume (and we do) that $Q'$ acts faithfully on $X'$. By $(1)$, by enlarging $\eta'$ if necessary, and by \cite[4.4.1]{Vi}, we can further assume (and we do) that there exists a surjective homomorphism $\varphi' : \tilde{Q}' \rightarrow Q'$ of finite groups and a proper surjective $\tilde{Q}'$-equivariant generically finite flat morphism $a : Y' \rightarrow X'$ with $Y'$ strictly semistable over $S'$, such that $k(Y')^{\Gamma}$ is purely inseparable over $k(X')$, where $\Gamma$ is the kernel of $\varphi'$. The profinite group $\tilde{Q} = \tilde{Q}' \times_{\Gal(\eta'/\eta_0)} \Gal(\eta/\eta_0)$ acts admissibly on $Y = Y'_{S}$, and is endowed with a natural surjective homomorphism $\varphi : \tilde{Q} \rightarrow Q$ with kernel $\Gamma$ such that $a_{S} : Y \rightarrow X$ is $\tilde{Q}$-equivariant.

 Let $j : U \rightarrow X'$ be a $Q'$-invariant dense open subscheme, with complement $i : Z \rightarrow X'$, such that $V = U \times_{X'} Y$ is finite flat over $U$ and such that $V \rightarrow V/\Gamma$ is finite \'etale. Let $j' : V \rightarrow Y$ be the canonical open immersion, with complement $i': T \rightarrow Y$. By the induction hypothesis, the conclusion of Lemma \ref{lem3.3.2} holds for the $Q$-equivariant morphism $Z_S \rightarrow S$ and for the $\tilde{Q}$-equivariant morphism $T_S \rightarrow S$. Thus, if $R\Psi_{Y'/S'} E$ is in $M^{\mathrm{eff}}(X_{s},Q,E)$, then so are $R\Psi_{Y'/S'} (j_{\eta'})_!E$ and $R\Gamma^I(R(a_s)_* R\Psi_{Y'/S'} (j_{\eta'})_!E) \simeq R\Psi_{X'/S'} (j_{\eta'})_!E$, in which case $R\Psi_{X'/S'}E$ is in $M^{\mathrm{eff}}(X_{s},Q,E)$ as well.

\item By $(4)$, we can assume (and we do) that $X'$ is strictly semistable over $S'$, purely of some relative dimension $d$. Since the problem is Zariski-local on $X'$, we can assume (and we do) that $X'$ is \'etale over the spectrum of $\Ow_{S'}[(Y_j)_{j \in R}]/(\prod_{j \in J} Y_j - \pi)$, for some uniformizer $\pi$ of $\Ow_{S'}$ and some finite sets $J \subseteq R$. Let $D_j$ be the closed subscheme of $X$ defined by the vanishing of $(\pi,Y_j)$, and let $Z^{(r)}$ be the disjoint union of $( \cap_{j \in J'} D_j)_{J' \subseteq J, |J'| = r +1}$. The finite group $Q'$ acts admissibly on $Z^{(r)}$, and the natural morphism $a^{(r)} : Z^{(r)} \rightarrow X'$ is $Q'$-equivariant. By the description of the monodromy filtration for the nearby cycles of strictly semistable $S'$-schemes, cf. \cite[2.7]{Sa03}, the object $R\Psi_{X'/S'} E$ of $D_{c}^b(X_{s},Q,\widehat{E})$ belongs to the sub-triangulated category generated by the objects
$$
 \bigoplus_{\substack{r \geq |t| \\ r = t \pmod 2}} Ra^{(r)}_* E\left(\frac{t-r}{2} \right)[-r].
$$
where $t$ is an integer. Since $M^{\mathrm{eff}}(X_{s},Q,E)$ is stable by shifts and by negative Tate twists, we obtain that $R\Psi_{X'/S'} E$ belongs to $M^{\mathrm{eff}}(X_{s},Q,E)$, hence concluding the proof of Lemma \ref{lem3.3.2}.

\end{enumerate}

\subsection{\label{3.20}} Let $S$ be a henselian trait, with closed point $s$ and generic point $\eta$, endowed with an admissible right action of a profinite group $Q$. Let $\overline{s}$ be a separable closure of $s$ with perfection $\overline{s}^{\mathrm{perf}}$ and let $\overline{\eta}$ be a separable closure of $\eta_{\overline{s}}$. Let us consider $S_0 = S/Q$, a trait with closed point $s_0$ and generic point $\eta_0$. Let $\overline{S}$ be the normalization of $S$ in $\overline{\eta}$ and let $\overline{Q}$ be the fiber product $\Gal(\overline{\eta}/\eta_0) \times_{\Gal(\eta/\eta_0)} Q$. Then $\overline{S}$ is a $\overline{Q}$-trait, with $\overline{S}/\overline{Q} = S_0$ and with closed point $\overline{s}^{\mathrm{perf}}$.

Let $X$ be an $S$-scheme of finite type. Then the nearby cycles functor
$$
R\Psi_{X/S} : D_c^b(X_{\eta},Q,\widehat{E}) \rightarrow D_c^b(X_{\overline{s}},\overline{Q},\widehat{E}),
$$
coincides with the composition
$$
D_c^b(X_{\eta},Q,\widehat{E}) \simeq D_c^b(X_{\overline{\eta}},\overline{Q},\widehat{E}) \xrightarrow[]{R\Psi_{X_{\overline{S}}/\overline{S}}} D_c^b(X_{\overline{s}^{\mathrm{perf}}},\overline{Q},\widehat{E}) \simeq D_c^b(X_{\overline{s}},\overline{Q},\widehat{E}) ,
$$
hence sends $M^{\mathrm{eff}}(X_{\eta},Q,E)$ to $M^{\mathrm{eff}}(X_{\overline{s}},\overline{Q},E)$ by Theorem \ref{teo3.3.1} and Proposition \ref{prop3.23.8}. Similarly, the nearby cycles functor
$$
R \Phi_{X/S} :  D_c^b(X,Q,\widehat{E}) \rightarrow D_c^b(X_{\overline{s}},\overline{Q},\widehat{E}),
$$
sends $M^{\mathrm{eff}}(X,Q,E)$ to $M^{\mathrm{eff}}(X_{\overline{s}},\overline{Q},E)$.

\begin{prop}\label{prop3.20.1} Let $S$ be a henselian trait, with closed point $s$ and generic point $\eta$, endowed with an admissible right action of a profinite group $Q$, such that $\eta$ is unramified over $\eta_0$. Let $X$ be a $Q$-equivariant $S$-scheme of finite type and let $j : X_{\eta} \rightarrow X$ be the canonical open immersion. Then the functor $Rj_*$ from $D_c^b(X_{\eta},Q,\widehat{E}) $ to $D_c^b(X,Q,\widehat{E})$ sends $M^{\mathrm{eff}}(X_{\eta},Q,E)$ to $M^{\mathrm{eff}}(X,Q,E)$.
\end{prop}

We keep the notation from \ref{3.20}. Let $i : X_s \rightarrow X$ be the canonical closed immersion. For any object $\F$ of $M^{\mathrm{eff}}(X_{\eta},Q,E)$, we have a distinguished triangle
$$
j_! \F \rightarrow Rj_* \F \rightarrow i_* i^*Rj_* \F  \xrightarrow[]{[1]},
$$
and $j_! \F$ is in $M^{\mathrm{eff}}(X,Q,E)$. We now prove that $i^*Rj_* \F$ is in $M^{\mathrm{eff}}(X_s,Q,E)$. Let $I_{\eta}$ be the interia subgroup in $\Gal(\overline{\eta}/\eta)$. Then the quotient $\overline{Q}/I_{\eta}$ coincides with $\Gal(\overline{s}/s_0) \times_{\Gal(s/s_0)} Q$. Let us consider the composition
$$
 D_c^b(X_{\overline{s}},\overline{Q},\widehat{E}) \xrightarrow[]{R \Gamma^{I_{\eta}}} D_c^b(X_{\overline{s}},\overline{Q}/ I_{\eta},\widehat{E}) \simeq D_c^b(X_{s},Q,\widehat{E}),
$$
which we still denote by $R \Gamma^{I_{\eta}}$. Then we have
$$
i^*Rj_* \F \simeq R \Gamma^{I_{\eta}}( R\Psi_{X/S} \F).
$$
By \ref{3.20}, it is therefore sufficient to prove that $R \Gamma^{I_{\eta}}$ sends $M^{\mathrm{eff}}(X_{\overline{s}},\overline{Q},E) $ to $M^{\mathrm{eff}}(X_{s},Q,E)$. Let us consider a $\overline{Q}$-equivariant separated morphism $a : Z \rightarrow X_{\overline{s}}$ of finite presentation, and an object $\G$ of $D_c^b(Z,\overline{Q},\Ow_E)$. Let us prove that $R \Gamma^{I_{\eta}}(Ra_! \G)$ is in $M^{\mathrm{eff}}(X_{s},Q,E)$. There exists an open subgroup $I'$ of $I_{\eta}$, normal in $\overline{Q}$, which acts trivially on $(Z,\F)$. Then there exists a $Q$-equivariant extension $\eta' \rightarrow \eta$ such that $I_{\eta'} = I'$. We have
$$
R \Gamma^{I_{\eta}}(Ra_! \G) \simeq R \Gamma^{I_{\eta}/ I_{\eta'}} R \Gamma^{I_{\eta'}}(Ra_! \G),
$$
and $R \Gamma^{I_{\eta}/ I_{\eta'}} $ has the desired property by Proposition \ref{prop3.26.10}. By replacing $\eta$ with $\eta'$ if necessary, we can assume (and we henceforth do) that $I_{\eta}$ acts trivially on $(Z,\F)$.

Let $p$ be the characteristic exponent of $s$, and let $P$ be the unique pro-$p$-Sylow subgroup of $I_{\eta}$. Then for any integer $n \geq 1$, we have
$$
R \Gamma^{I_{\eta}}(Ra_! (\G \otimes_{\Ow_E}^L \Ow_E/n\Ow_E)) \simeq R \Gamma^{I_{\eta}/P}(Ra_! (\G \otimes_{\Ow_E}^L \Ow_E/n\Ow_E)),
$$
since both terms vanish unless $n$ is prime to $p$, in which case it follows from the fact that $\Gamma^P$ is an exact functor on prime to $p$ torsion coefficients. The profinite group $I_{\eta}/P$ is canonically isomorphic to the procyclic group $\prod_{\ell \neq p} \mathbb{Z}_{\ell}(1)$, where the Tate twist reflects the conjugation action of $\overline{Q}/I_{\eta}$. For any object $\mathcal{H}$ of $D_c^b(X_{\overline{s}}, \overline{Q}/I_{\eta}, \widehat{E})$, we have a canonical isomorphism
$$
R \Gamma^{I_{\eta}/P}(\mathcal{H})  \simeq \mathcal{H} \oplus \mathcal{H}(-1)[-1].
$$
Since $M^{\mathrm{eff}}(X_{s},Q,E)$ is stable by negative Tate twists, by shifts and by direct sums, we deduce that $R \Gamma^{I_{\eta}/P}(Ra_! \G)$ is in $M^{\mathrm{eff}}(X_{s},Q,E)$, hence the conclusion.

\section{Six functors formalism for schemes of finite type over a field  \label{full}}

Let $E$ be a number field, with ring of integers $\Ow_E$. Let $s$ be the spectrum of a field, endowed with an admissible right action from a profinite group $Q$.

\subsection{\label{6.3}} We start with the following straightforward application of Proposition \ref{prop3.20.1}.

\begin{prop}\label{prop6.3.1} Let $f : X \rightarrow Y$ be a $Q$-equivariant morphism of $s$-schemes of finite type. Then the functor $Rf_*$ from $D_c^b(X,Q,\widehat{E}) $ to $D_c^b(Y,Q,\widehat{E})$ sends $M^{\mathrm{eff}}(X,Q,E)$ to $M^{\mathrm{eff}}(Y,Q,E)$.
\end{prop}

We first consider the case where $f$ is an open immersion. By replacing $Y$ with a suitable $X$-admissible blow-up if necessary, we can assume (and we do) that the complement of $X$ in $Y$ is the support of a Cartier divisor $D$. Let $\F$ be an object of $M^{\mathrm{eff}}(X,Q,E)$ and let $U$ be a $Q$-invariant affine open subscheme of $Y$ such that the closed immersion $i : D \cap U \rightarrow U$ is defined by a single equation, i.e. is the fiber above $0$ of a morphism $g : U \rightarrow \mathbb{A}^1_s$. By applying Proposition \ref{prop3.20.1} to the pullback of $g$ along the henselization of $\mathbb{A}^1_s$ at $0$, we obtain that $i^*( Rf_* \F_{|U})$ belongs to $M^{\mathrm{eff}}(U \cap D,Q,E)$. Thus $Rf_* \F_{|U}$ belongs to $M^{\mathrm{eff}}(U,Q,E)$, and the conclusion follows.

We now consider the case where $f$ is separated. We can find a dense $Q$-equivariant open immersion $j : X \rightarrow \overline{X}$ with $\overline{f} : \overline{X} \rightarrow Y$ a $Q$-equivariant proper morphism. We have $Rf_* \simeq R\overline{f}_* Rj_*$ and the conclusion follows from the case of open immersions already handled, and the case of proper morphisms from Proposition \ref{prop3.23.2}.

Let us now prove the general case. We can assume (and we do) that $Y$ is affine. We write $X$ as the union of a finite collection $(U_r)_{r \in R}$ of quasi-affine open subschemes. We argue by induction on the cardinality of $R$, the case where $R$ is empty being trivial. If $R$ is non empty, let $r_0$ be an element of $R$, and let $U = U_{r_0}$ and $V = \cup_{r \neq r_0} U_r$, with corresponding open immersions $j_U : U \rightarrow X$ and $j_V: V \rightarrow X$. Then for an object $\F$ of $M^{\mathrm{eff}}(X,Q,E)$ we have a distinguished triangle
$$
Rf_* \F \rightarrow R(fj_U)_* \F \oplus R(fj_V)_* \F \rightarrow  R(fj_{U \cap V})_* \F \xrightarrow[]{[1]}.
$$
By the induction hypothesis, both $R(fj_V)_* \F$ and $R(fj_{U \cap V})_* \F$ are in $M^{\mathrm{eff}}(Y,Q,E)$. Moreover, the already obtained result for separated morphisms implies that $R(fj_U)_* \F$ belongs to $M^{\mathrm{eff}}(Y,Q,E)$. Thus $Rf_* \F$ is in $M^{\mathrm{eff}}(Y,Q,E)$.

\subsection{\label{6.1}} Let $X$ be a $Q$-scheme. The category $M(X,Q,E)$ of \textit{$Q$-equivariant $E$-systems} on $X$ is the full subcategory of $D_{c}^b(X,Q,\widehat{E})$ (cf. \ref{3.27}) consisting of $Q$-equivariant $\widehat{E}$-systems $\F$ on $X$ with the following property: for any point $x$ of $X$, there exists a $Q$-invariant affine open subscheme $U$ and an integer $d$ such that such that the Tate twist $\F(-d)_{|U}$ belongs to $M^{\mathrm{eff}}(U,Q,E)$. 

If $X$ is quasicompact, then $M(X,Q,E)$ simply consists of objects $\F$ of $D_{c}^b(X,Q,\widehat{E})$ such that $\F(-d)$ belongs to $M^{\mathrm{eff}}(X,Q,E)$ for some integer $d$. 

\subsection{\label{6.2}} For a $Q$-equivariant separated $s$-scheme $X$ of finite type, we denote by $\tau_X : X \rightarrow s$ the structural morphism, and by $D_{X/s} = R Hom( -  , R\tau_X^{!} E)$ the duality functor from $D_{c}^b(X,Q,\widehat{E})$ to itself. The functor $D_{X/s}$ is an equivalence of categories, with itself as a quasi-inverse.

\begin{prop}\label{prop6.2.1} Let $X$ be a $Q$-equivariant separated $s$-scheme of finite type. Then the duality functor $D_{X/s}$ from $D_c^b(X,Q,\widehat{E}) $ to itself preserves $M(X,Q,E)$.
\end{prop}

By Proposition \ref{prop3.19.1}, it is enough to consider the image by $D_{X/s}$ of objects of the form
$$
\F = (V \otimes_E Ra_! E)^I,
$$
where $a : Z \rightarrow X$ is a $Q'$-equivariant separated morphism of finite type, for some extension $Q'$ of $Q$ by a finite group $I$, and where $V$ is a finite dimensional continuous $E$-representation of $Q'$. 

We prove by induction on the dimension $d$ of $Z$ that $D_{X/s}(\F)$ is in $M(X,Q,E)$. We can assume that $Z$ is reduced. By replacing $s$ with a suitable purely inseparable extension, we may further assume (and we do) that $X$ is geometrically reduced. Then there exists a dense open immersion $j : U \rightarrow Z$ such that $U$ is smooth over $s$. If $i :T \rightarrow Z$ is a closed immersion with support given by the complement of $U$ in $Z$, then we have a distinguished triangle
$$
R\Gamma^I(V \otimes_E R(aj)_! E) \rightarrow \F \rightarrow R\Gamma^I(V \otimes_E R(ai)_! E)  \xrightarrow[]{[1]}.
$$
By the induction hypothesis, the functor $D_{X/s}$ sends $R\Gamma^I(V \otimes_E R(ai)_! E) $ to an object of $M(X,Q,E)$. Thus, by replacing $a$ with $aj$ if necessary, we can assume (and we do) that $Z$ is smooth over $s$. We can further assume (and we do) that $Z$ is connected, purely of dimension $d$. We then have
$$
D_{X/s} \F \simeq (V^{\vee} \otimes_E D_{X/s}  Ra_! E)^I \simeq (V^{\vee} \otimes_E Ra_* D_{Z/s}   E)^I \simeq (V^{\vee} \otimes_E Ra_*  E(d)[2d])^I,
$$
and the conclusion follows from Proposition \ref{prop6.3.1} applied to $a$.

\begin{cor}\label{cor6.2.2} Let $f : X \rightarrow Y$ be a $Q$-equivariant morphism of separated $s$-schemes of finite type. Then the functor $Rf^!$ from $D_c^b(Y,Q,\widehat{E}) $ to $D_c^b(X,Q,\widehat{E})$ sends $M(X,Q,E)$ to $M(Y,Q,E)$.
\end{cor}

This follows from Proposition \ref{prop6.2.1}, since we have $Rf^! \simeq D_{X/s} f^* D_{Y/s}$.

\begin{cor}\label{cor6.2.3} Let $X$ be a $Q$-equivariant $s$-scheme of finite type. Then the functor $RHom$ from $D_c^b(X,Q,\widehat{E}) \times D_c^b(X,Q,\widehat{E})$ to $D_c^b(X,Q,\widehat{E})$ sends $M(X,Q,E) \times M(X,Q,E)$ to $M(X,Q,E)$.
\end{cor}

Since the assertion is Zariski-local on $X$, we can assume (and we do) that $X$ is separated over $s$, in which case it follows from Propositions \ref{prop6.2.1} and \ref{prop3.19.10}, since we have 
$$
RHom(\F,\G) \simeq D_{X/s}( \F \otimes D_{X/s}(\G) ).
$$

\section{Linearized $\varepsilon$-factors \label{vareps}}

Let $E$ be a number field, with ring of integers $\Ow_E$. Let $S$ be a henselian trait of positive characteristic $p$, with perfect closed point $s$ and generic point $\eta$, endowed with an admissible right action from a profinite group $Q$. Let $\overline{s}$ be a separable closure of $s$ and let $\overline{\eta}$ be a separable closure of $\eta_{\overline{s}}$.

We assume throughout this section that $S$ is unramified over $S_0 = S/Q$ and that $S$ is an $s$-scheme. We denote by $s_0$ and $\eta_0$ the closed point and generic point of $S_0$ respectively. Let $\overline{S}$ be the normalization of $S$ in $\overline{\eta}$ and let $\overline{Q}$ be the fiber product $\Gal(\overline{\eta}/\eta_0) \times_{\Gal(\eta/\eta_0)} Q$. Then $\overline{S}$ is a $\overline{Q}$-trait, with $\overline{S}/\overline{Q} = S_0$ and with closed point $\overline{s}^{\mathrm{perf}}$. We denote by $I_{\eta}$ the inertia group $\Gal(\overline{\eta}/\eta_{\overline{s}})$ of $\eta$.

\subsection{\label{3.21.1}} Let $\pi$ be a uniformizer on $S_0$. Then we have a homomorphism from $k(s)[T,T^{-1}]$ to $k(\eta)$ sending $T$ to $\pi$, and we still denote  by $\pi$ the corresponding $Q$-equivariant morphism
$$
\pi : \eta \rightarrow \mathbb{G}_{m,s}.
$$
The theory of Gabber-Katz extensions, cf. \cite[4.1]{G192}, provides a fully faithful functor
$$
\mathrm{GK}_{\pi} : \Sh(\eta, Q,  \Lambda) \rightarrow \Sh(\mathbb{G}_{m,s}, Q,   \Lambda),
$$
for any ring $\Lambda$, such that the composition $\pi^*\mathrm{GK}_{\pi}$ is naturally isomorphic to the identity functor. For any $Q$-equivariant $\Lambda$-sheaf $\F$ on $\eta$, the extension $\mathrm{GK}_{\pi}(\F)$ is a local system tamely ramified at $\infty$, and its geometric monodromy group admits a unique $p$-Sylow subgroup. These properties characterize the essential image of $\mathrm{GK}_{\pi}$.

\subsection{\label{3.21.2}} Let $\psi : \mathbb{F}_p \rightarrow E^{\times}$ be a non trivial homomorphism, and let $\Lc_{\psi}$ be the corresponding Artin-Schreier $E$-local system on the affine line $\mathbb{A}_{\mathbb{F}_p}^1$. For any morphism $f$ from a scheme to $\mathbb{A}_{
\mathbb{F}_p}^1$, we denote by $\Lc_{\psi} \{ f \}$ the pullback $f^*\Lc_{\psi}$. In \cite[4.13]{G192}, we considered the functor
\begin{align*}
\Art_{\pi, \psi} :\Sh(\overline{s}, \overline{Q},  \widehat{\Ow}_E) = \Sh(\eta, Q,  \widehat{\Ow}_E) &\rightarrow  \Sh(\overline{s}, \overline{Q}/I_{\eta},   \widehat{\Ow}_E) = \Sh(s, Q,  \widehat{\Ow}_E)\\
\F &\mapsto R\Gamma_c( \mathbb{G}_{m,\overline{s}}, \Lc_{\psi} \{ -t \} \otimes_E \mathrm{GK}_{\pi}(\F) )[1],
\end{align*}
where $t$ is the coordinate on $\mathbb{G}_{m,s}$, and the complex $\Art_{\pi,\psi}(\F)$ is concentrated in degree $0$. Moreover, it is shown in \cite[2.18]{G192} how to construct a functor
$$
\Art_{\pi, \psi} : \Sh(X_{\overline{s}}, \overline{Q},   \widehat{\Ow}_E) \rightarrow \Sh(X_{\overline{s}}, \overline{Q}/I_{\eta},   \widehat{\Ow}_E) \simeq \Sh(X, Q,   \widehat{\Ow}_E),
$$
for any $Q$-equivariant $s$-scheme $X$, which coincides with the above functor for $X=s$, and is compatible with pullbacks, pushforwards, and their derived functors.

\subsection{\label{3.21.4}} In the case $X=s$, the functor $\Art_{\pi, \psi}$ can be considered as a linearization of the local $\varepsilon$-factor. Indeed, for any $\F$ in $\Sh(\eta, 1, \widehat{E})$, and any non archimedean place $\lambda$ of $E$ not dividing $p$, the $\lambda$-adic component $\Art_{\pi, \psi}(\F)_{\lambda}$ of $\Art_{\pi, \psi}(\F)$ satisfies by \cite[4.14]{G192} the identities
\begin{align}\label{idento}
\rk(\Art_{\pi, \psi}(\F)_{\lambda}) &= a(S,j_! \F_{\lambda}) \\
\det(\Art_{\pi, \psi}(\F)_{\lambda}) &= \varepsilon_{\psi}(S,j_!\F_{\lambda}, d\pi)
\end{align}
where $j : \eta \rightarrow S$ is the canonical open immersion, where $a(S, -)$ is the Artin conductor or total dimension, and $\varepsilon_{\psi}(S, -,d\pi)$ is the local $\varepsilon$-factor from \cite{G19}.

\subsection{\label{3.21.5}} Let $I'$ be an open subgroup of $I_{\eta}$, which is normal in $\overline{Q}$. We have an exact sequence
$$
1 \rightarrow \Gamma \rightarrow \overline{Q}/ I' \rightarrow Q \rightarrow 1,
$$
where $\Gamma = \Gal(\overline{\eta}/\eta)/I'$. Let $\eta' \rightarrow \eta_{\overline{s}}$ be a finite separable $\overline{Q}/ I'$-equivariant extension, which is Galois of group $I_{\eta}/I'$. Let us consider its Gabber-Katz extension $a : Z' \rightarrow \mathbb{G}_{m,\overline{s}}$, which fits into a $\overline{Q}/ I'$-equivariant diagram
\begin{center}
 \begin{tikzpicture}[scale=1]

\node (A) at (0,0) {$Z'$};
\node (C) at (2,0) {$\mathbb{G}_{m,\overline{s}}.$};
\node (D) at (0,1) {$ \eta'$};
\node (E) at (2,1) {$ \eta$};

\path[->,font=\scriptsize]
(A) edge  (C)
(D) edge (A)
(E) edge (C)
(D) edge (E);
\end{tikzpicture} 
\end{center}

Then for any $Q$-equivariant $s$-scheme $X$ and any element $\F$ of $\Sh(X_{\overline{s}}, \overline{Q}/I',   \widehat{E})$, we have a functorial isomorphism
\begin{align*}
\Art_{\pi, \psi}(\F) &\simeq  \left( R\Gamma_c( \mathbb{G}_{m,\overline{s}}, \Lc_{\psi} \{ - t \} \otimes_E a_* E )[1] \otimes_{E} \F \right)^{I_{\eta}/I'} \\
&\simeq \left( R\Gamma_c(Z', \Lc_{\psi} \{ - a \} )[1] \otimes_{E} \F \right)^{I_{\eta}/I'}.
\end{align*}
For an object $\F$ of $D_c^b(X_{\overline{s}}, \overline{Q}/I', \widehat{E})$, we obtain an isomorphism
$$
\Art_{\pi, \psi}(\F) = \left( \Mac_{I'}\otimes_{E} \F \right)^{I_{\eta}/I'},
$$
in $D_c^b(X_{\overline{s}}, \overline{Q}/I_{\eta},   \widehat{E})$, where $\Mac_{I'}$ is the element of $M^{\mathrm{eff}}(X_{\overline{s}}, \overline{Q}/I',   E)$ defined by $ \tau_{X_{\overline{s}}}^*(\tau_{Z'})_!  \Lc_{\psi} \{ - a \}[1]$, where $\tau_{Z'}, \tau_{X_{\overline{s}}}$ are the structural morphisms from $Z'$ and $X_{\overline{s}}$ to $\overline{s}$.

\begin{prop}\label{propvareps} Let $X$ be a $Q$-equivariant $s$-scheme and let $\psi : \mathbb{F}_p \rightarrow E^{\times}$ be a non trivial homomorphism. Then the functor $\Art_{\pi, \psi}$ from $D_c^b(X_{\overline{s}}, \overline{Q},   \widehat{E})$ to $D_c^b(X_s,Q,\widehat{E})$ sends $M^{\mathrm{eff}}(X_{\overline{s}}, \overline{Q},E)$ to $M^{\mathrm{eff}}(X, Q,E)$.
\end{prop}

Indeed, by Proposition \ref{prop3.16.9}, it is enough to prove that for any open subgroup $I'$ of $I_{\eta}$ which is normal in $\overline{Q}$, the functor $\Art_{\pi, \psi}$ sends objects of $M^{\mathrm{eff}}(X_{\overline{s}}, \overline{Q}/I',E)$ to $M^{\mathrm{eff}}(X, Q,E)$, and the conclusion then follows from \ref{3.21.5}.

\section{Independence of $\ell$ \label{indep}}

Let $E$ be a number field, with ring of integers $\Ow_E$. We denote by $\mathbb{L}$ the set of non-archimedean places of $E$. For any $\lambda$ in $\mathbb{L}$, we denote by $E_{\lambda}$ the $\lambda$-adic completion of $E$ and by $\Ow_{E_{\lambda}}$ its ring of integers. 

For a $Q$-scheme $X$, an object $\F$ of $D_c^b(X,Q,\widehat{E})$, and an element $\lambda$ of $\mathbb{L}$ of residue characteristic $\ell_{\lambda}$, we denote by $\F_{\lambda}$ the $Q$-equivariant $\lambda$-adic complex on $X\left[ \ell_{\lambda}^{-1} \right]$ associated to $\F$.

\subsection{\label{5.8}} The following result, which asserts that $E$-systems satisfy the classical notion of $E$-compatibility, is a refinement of \cite[3.3]{DL76}, but the proof is mostly identical. We therefore do not claim any originality here. 

\begin{prop}\label{invar} Let $Q$ be a profinite group. Let $s$ be the spectrum of a finite field of characteristic $p$, endowed with the trivial action of $Q$. Let $\overline{s}$ be a separable closure of $s$, and let $\Frob_s$ be the geometric Frobenius in $\Gal(\overline{s}/s)$. Let $\F = (\F_n)_n$ be an object of $M(s,Q,E)$ (cf. \ref{6.1}).
\begin{enumerate}
\item Let $q$ be an element of $Q$. Then there exists a unique element $\Tr(q| \F_{\overline{s}})$ of $\Ow_E \left[ p^{-1} \right]$ such that for any element $\lambda$ of  $\mathbb{L}$ not dividing $p$, the trace
$$
\Tr(q |  \F_{\lambda, \overline{s}}),
$$
coincides with the image of $\Tr(q| \F_{\overline{s}})$ in $\Ow_{E_{\lambda}}$.
\item There exists a unique rational function $ \det(1 - T \Frob_s | \F_{ \overline{s}})$ in $E(T)$ such that for any $\lambda$ in $\mathbb{L}$, we have
$$
r(T) = \det(1 - T \Frob_s | \F_{\lambda, \overline{s}}),
$$
in $E_{\lambda}(T)$. Moreover, the rational function $ \det(1 - T \Frob_s | \F_{ \overline{s}})$ is a quotient of polynomials in $1 + T \Ow_E \left[ p^{-1} \right][T]$ with invertible leading coefficients.

\end{enumerate}
\end{prop}

We follow closely the proof of \cite[3.3]{DL76}. Since these properties are stable by Tate twist, we can assume that $\F$ is effective. We can further assume (and we do) that $\F$ is of the form $Ra_! \F$, for some $Q$-equivariant separated morphism $a : Z \rightarrow s$ of finite type and some object $\G$ of $D_c^b(Z,Q,)$. For each $\lambda$ in $\mathbb{L}$, we fix an algebraic closure $\overline{E}_{\lambda}$ of $E_{\lambda}$. For any element $c$ of $\overline{E}_{\lambda}$, let us set
$$
\alpha_{c,\lambda}(q) = \Tr(q | R\Gamma(Z_{\overline{s}}, \G \otimes_E \overline{E}_{\lambda})^{\Frob_s \sim c}),
$$
where $\sim c$ indicates the generalized eigenspace for the eigenvalue $c$. For each integer $n \geq 1$, we have 
$$
\sum_{c \in \overline{E}_{\lambda}^{\times}} \alpha_{c,\lambda}(q) c^n = \Tr( q \Frob_s^n  | R\Gamma(Z_{\overline{s}},\G \otimes_E \overline{E}_{\lambda})).
$$
If $s_n$ is the unique extension of degree $n$ of $s$ in $\overline{s}$, then the automorphism $q \Frob_s^n $ of $Z_{\overline{s}}$ is the geometric Frobenius automorphism of some $s_n$-form of $Z_{\overline{s}}$, and the Grothendieck-Lefschetz trace formula therefore gives
$$
\Tr( q \Frob_s^n  | R\Gamma(Z_{\overline{s}},\G \otimes_E \overline{E}_{\lambda}))  = \sum_{z \in Z(\overline{s})^{q \Frob_s^n  = 1}} \Tr( \Frob_z | \G_{\overline{z}}).
$$
In particular, this an element of $\Ow_E$. Thus, for any $E$-linear field isomorphism $\iota : \overline{E}_{\lambda} \rightarrow \overline{E}_{\lambda'}$, we have
$$
\sum_{c \in \overline{E}_{\lambda'}^{\times}} \alpha_{c,\lambda'}(q) c^n = \iota \left( \sum_{c \in \overline{E}_{\lambda}^{\times}} \alpha_{c,\lambda}(q) c^n \right) = \sum_{c \in \overline{E}_{\lambda}^{\times}} \iota(\alpha_{c,\lambda}(q)) \iota(c)^n,
$$
and consequently $\iota(\alpha_{c,\lambda}(q)) = \alpha_{\iota(c),\lambda'}(q)$ by linear independence of the maps $(n \mapsto c^n)_{c \in \overline{E}_{\lambda'}^{\times}}$. Consequently, the sum
$$
\Tr(q | R\Gamma(Z_{\overline{s}},\G \otimes_E \overline{E}_{\lambda}) ) = \sum_{c \in \overline{E}_{\lambda}^{\times}} \alpha_{c,\lambda}(q),
$$
is invariant by $\Aut(\overline{E}_{\lambda}/E)$, hence belongs to $E(T)$, and is independent of $\lambda$. Similarly, we have
$$
\iota \det(1 - T \Frob_s | \F_{\lambda, \overline{s}}) = \iota \prod_{c \in \overline{E}_{\lambda}^{\times}} (1 - T c)^{\alpha_{c,\lambda}(1)} = \det(1 - T \Frob_s | \F_{\lambda', \overline{s}}).
$$
Thus the rational function $\det(1 - T \Frob_s | \F_{\lambda, \overline{s}})$ is $\Aut(\overline{E}_{\lambda}/E)$-invariant, hence belongs to $E(T)$, and is independent of $\lambda$ in $\mathbb{L}$. 

\begin{rema} The proof yields the more general statement that the element
$$
\prod_{c \in \overline{E}_{\lambda}^{\times}} (1 - T c)^{\alpha_{c,\lambda}(q)},
$$
of $\overline{E}_{\lambda}(T)^{\times} \otimes_{\mathbb{Z}} \mathbb{Z}[\mu_{\infty}(\overline{E}_{\lambda})]$ is actually in
$$
(\overline{E}(T)^{\times}  \otimes_{\mathbb{Z}} \mathbb{Z}[\mu_{\infty}(\overline{E})])^{\Gal(\overline{E}/E)},
$$
and is independent of $\lambda$.
\end{rema}

\begin{cor}\label{invar3} Let $Q$ be a profinite group. Let $s$ be the spectrum of a field, endowed with the trivial action of $Q$, with separable closure $\overline{s}$. Let $\F$ be an object of $M(s,Q,E)$ (cf. \ref{6.1}). Then for any element $q$ of $Q$, there exists a unique element $\Tr(q| \F_{\overline{s}})$ of $\Ow_E \left[ p^{-1} \right]$ such that for any element $\lambda$ of $\mathbb{L}$ of residual characteristic invertible on $s$, the trace
$$
\Tr(q |  \F_{\lambda, \overline{s}}),
$$
coincides with the image of $\Tr(q| \F_{\overline{s}})$ in $\Ow_{E_{\lambda}}$.
\end{cor}

Indeed, there exists a dominant morphism $s \rightarrow S$ to an integral scheme $S$ of finite type over $\mathbb{Z}$, such that $\F$ is the pullback on object $\F'$ of $M(S,Q,E)$. By Proposition \ref{prop3.23.3} and by shrinking $S$ if necessary, we can assume (and we do) that $\F'_{\lambda}$ is a lisse $\lambda$-adic complex for any $\lambda$. The result then follows by applying Proposition \ref{invar}(1) to a closed point of $S$.

\begin{cor}\label{invarcharccycle} Let $s$ be the spectrum of a perfect field of characteristic exponent $p \geq 1$, let $X$ be a smooth $s$-scheme of pure dimension $d$, and let $\F$ be an object of $M(X,1,E)$. Then the characteristic cycle $\mathrm{CC}(\F_{\lambda})$ from \cite[Th. 4.9]{Sai} (cf. \ref{intro5}) does not depend on the choice of an element $\lambda$ of $\mathbb{L}$ not dividing $p$.
\end{cor}

Let $\F$ be an object of $M(X,1,\Ow_E)$. Let us consider a commutative diagram
\begin{center}
 \begin{tikzpicture}[scale=1]

\node (A) at (0,0) {$S$};
\node (C) at (2,0) {$s$};
\node (D) at (0,1) {$U$};
\node (E) at (2,1) {$X$};

\path[->,font=\scriptsize]
(A) edge  (C)
(D) edge node[left]{$f$} (A)
(E) edge (C)
(D) edge node[above]{$j$} (E);
\end{tikzpicture} 
\end{center}
where $j$ is \'etale and $S$ is a smooth curve over $s$. Let $u$ be a closed point of $U$, with separable closure $\overline{u}$, and let $\eta$ be the generic point of the strictly henselian trait $S_{(f(\overline{u}))}$, with separable closure $\overline{\eta}$. Then the stalk of $R \Phi_{X \times_S S_{(f(\overline{u}))}/ S_{(f(\overline{u}))}} j^* \F$ at $\overline{u}$ is in $M(\overline{u}, \Gal(\overline{\eta}/\eta), E)$ (cf. \ref{3.20}). By Proposition \ref{invar3}, the trace
$$
\Tr(\sigma | (R \Phi_{X \times_S S_{(f(\overline{u}))}/ S_{(f(\overline{u}))}} j^* \F_{\lambda})_{\overline{u}}),
$$
is, for each $\sigma$ in the inertia group $\Gal(\overline{\eta}/\eta)$, an element of $E$ which is independent of $\lambda$ in $\mathbb{L}$ not dividing $p$. In particular, the Artin conductor (or total dimension) of $(R \Phi_{X \times_S S_{(f(\overline{u}))}/ S_{(f(\overline{u}))}} j^* \F_{\lambda})_{\overline{u}}$ does not depend on $\lambda$. For any two elements $\lambda,\lambda'$ of $\mathbb{L}$ not dividing $p$, there exists a closed conical subset $C$ of $T^*X$ of pure dimension $d$ on which both $\F_{\lambda}$ and $\F_{\lambda'}$ are micro-supported. Then we can apply the preceding discussion to the particular situation where $u$ is an at most $C$-isolated characteristic point of $f$, cf. \cite[4.3]{Sai}, in which case we obtain  an equality
$$
(\mathrm{CC}(\F_{\lambda}),df)_{T^*U,u} = (\mathrm{CC}(\F_{\lambda'}),df)_{T^*U,u},
$$
of intersection numbers. By the uniqueness statement in \cite[Th. 4.9]{Sai}, we conclude that $\mathrm{CC}(\F_{\lambda})$ and $\mathrm{CC}(\F_{\lambda'})$ coincide.

\subsection{\label{5.10}} Let $s$ be the spectrum of a field of positive characteristic $p$, with separable closure $\overline{s}$. We define a topological group $W(\overline{s}/s)$, the Weil group of $\overline{s}/s$, as the fiber product
$$
W(\overline{s}/s) = W(\overline{\mathbb{F}}_p/\mathbb{F}_p) \times_{\Gal(\overline{\mathbb{F}}_p/\mathbb{F}_p)} \Gal(\overline{s}/s).
$$
of topological groups, where $\overline{\mathbb{F}}_p$ is the algebraic closure of $\mathbb{F}_p$ in $k(\overline{s})$ and $ W(\overline{\mathbb{F}}_p/\mathbb{F}_p)$ is the subgroup of $\Gal(\overline{\mathbb{F}}_p/\mathbb{F}_p)$ generated by the Frobenius element $\Frob_p$, endowed with the discrete topology. If the algebraic closure of $\mathbb{F}_p$ in $k(s)$ is a finite field, then the natural continuous homomorphism
$$
W(\overline{s}/s)  \rightarrow \Gal(\overline{s}/s),
$$
has dense image.

\subsection{\label{5.9}} Let $p \geq 1$ be an integer. By Dirichlet's unit theorem, the abelian group $\Ow_E[p^{-1}]^{\times}$ is finitely generated. Let $\widehat{\Ow_E[p^{-1}]^{\times}}$ be its profinite completion. For any finite set $\Sigma$ of non archimedean places of $E$ containing the places dividing $p$, the natural homomorphism
$$
\widehat{\Ow_E[p^{-1}]^{\times}} \rightarrow \prod_{\lambda \in \mathbb{L}  \setminus \Sigma} \Ow_{E_{\lambda}}^{\times},
$$
is injective by the weak Grunwald-Wang theorem \cite[IX.1 Th. 1]{AT}. 

\begin{prop}\label{invar2} Let $s$ be the spectrum of a field of positive characteristic $p$, with separable closure $\overline{s}$. Let $\F$ be an object of $M(s,1,E)$.

\begin{enumerate}
\item There exists a unique continuous homomorphism 
$$
\det(\F_{\overline{s}}) : \Gal(\overline{s}/s) \rightarrow  \widehat{\Ow_E[p^{-1}]}^{\times},
$$
such that for any $\lambda$ in $\mathbb{L}$ not dividing $p$, the composition
$$
\Gal(\overline{s}/s) \xrightarrow[]{\det(\F_{\overline{s}})}  \widehat{\Ow_E[p^{-1}]}^{\times} \rightarrow \Ow_{E_{\lambda}}^{\times},
$$
coincides with $\det(\F_{\lambda,\overline{s}})$. 
\item If the algebraic closure of $\mathbb{F}_p$ in $k(s)$ is a finite field, then the restriction of $\det(\F_{\overline{s}})$ to the Weil group $W(\overline{s}/s)$ uniquely factors through a continuous homomorphism $W(\overline{s}/s) \rightarrow  \Ow_E[p^{-1}]^{\times}$, where $ \Ow_E[p^{-1}]^{\times}$ is endowed with the discrete topology.
\end{enumerate}
\end{prop}

By Proposition \ref{prop3.23.5}, we can assume (and we do) that $k(s)$ is finitely generated over $\mathbb{F}_p$, in which case the algebraic closure $k$ of $\mathbb{F}_p$ in $s$ is a finite field. By Proposition \ref{prop3.23.5} again, we can consider a dominant morphism $s \rightarrow S$ to a normal geometrically connected $k$-scheme of finite type $S$, such that $\F$ is the pullback on object $\F'$ of $M(S,Q,E)$. By Proposition \ref{prop3.23.3} and by shrinking $S$ if necessary, we can assume (and we do) that $\F'_{\lambda}$ is a lisse $\lambda$-adic complex for any $\lambda$. Then for any $\lambda$ of residual characteristic $\ell$ prime to $p$, the homomorphism $\det(\F_{\lambda,\overline{s}})$ factors as
$$
\Gal(\overline{s}/s)^{\mathrm{ab}} \rightarrow \pi_1(S , \overline{s})^{\mathrm{ab}} \xrightarrow[]{\det(\F_{\lambda,\overline{s}}')} \Ow_{E_{\lambda}}^{\times}.
$$
By the Chebotarev density theorem, the classes $(\Frob_x)_x$, where $x$ runs over closed points of $S$, form a dense subset of $\pi_1(S , \overline{s})^{\mathrm{ab}}$. In particular, their images $(\Frob_k^{[k(x):k]})_x$ in $\Gal(\overline{k}/k)$ form a dense subset of $\Gal(\overline{k}/k)$. Thus the collection $([k(x):k])_x$ of integers is relatively coprime, and we can write
$$
1 = \sum_{i} n_i [k(x_i): k],
$$
for some relative integers $(n_i)_i$ and closed points $(x_i)_{i}$ of $S$. Let us define
$$
c = \prod_i \det(\Frob_{x_i} | \F'_{\overline{x}_i})^{n_i},
$$
which is an element of $\Ow_E[p^{-1}]^{\times}$, cf. \ref{invar}. By \cite[1.3.1]{De80}, for each $\lambda$ in $\Lambda$, we have
$$
\det(\F'_{\lambda,\overline{s}}) = \chi_{\lambda} c^{\mathrm{deg}},
$$
where $\mathrm{deg}$ is the composition of the homomorphism $\pi_1(S,\overline{s}) \rightarrow \Gal(\overline{k}/ k)$ with the isomorphism $\Gal(\overline{k}/ k) \simeq \widehat{\mathbb{Z}}$ which sends $\Frob_{k}$ to $1$, and where $\chi_{\lambda} : \pi_1(S,\overline{s}) \rightarrow \Ow_{E_\lambda}^{\times}$ is a character of finite order. For each closed point $x$ of $S$, we have
$$
\chi_{\lambda}(\Frob_x) =  \det(\Frob_{x} | \mathcal{N}_{\overline{x}})c^{-[k(x):k]},
$$
hence $\chi_{\lambda}(\Frob_x)$ belongs to $\Ow_E[p^{-1}]^{\times}$ and is independent of $\lambda$. By the aforementioned Chebotarev density theorem, this implies that the character $\chi_{\lambda}$ take its values in $\Ow_E[p^{-1}]^{\times}$ and does not depend of $\lambda$; let us denote it by $\chi$. The composition
$$
\Gal(\overline{s}/s)^{\mathrm{ab}} \rightarrow \pi_1(S , \overline{s})^{\mathrm{ab}} \xrightarrow[]{\chi c^{\mathrm{deg}} } \widehat{\Ow_E[p^{-1}]}^{\times},
$$
then satisfies the conclusion of $(1)$, and its restriction to $W(\overline{s}/s)$ factors through the continuous homomorphism $\chi c^{\mathrm{deg}}  : W(\overline{s}/s) \rightarrow \Ow_E[p^{-1}]^{\times}$.

\begin{cor}\label{invarepscycle} Let $s$ be the spectrum of a perfect field of positive characteristic $p$, with separable closure $\overline{s}$. We assume that $E$ contains a non trivial $p$-th root of unity. Let $X$ be a smooth $s$-scheme of pure dimension $d$, and let $\F$ be an object of $M(X,1,E)$. Then there exists a unique $d$-cycle $\mathcal{E}(\F)$ on the cotangent scheme $T^*X$, with coefficients in 
$$\mathbb{Q} \otimes_{\Z} \Hom_{\mathrm{cont}}(\Gal(\overline{s}/s), \widehat{\Ow_E[p^{-1}]}^{\times}/\mu_E),$$
such that for each element $\lambda$ of $\mathbb{L}$ not dividing $p$, the $\varepsilon$-cycle $\mathcal{E}(\F_{\lambda})$ from \cite[Th. 5.4]{Ta} (cf. \ref{intro7}) is obtained from $\mathcal{E}(\F)$ through the homomorphism $ \widehat{\Ow_E[p^{-1}]}^{\times}/\mu_E \rightarrow \widehat{\Ow_{E_{\lambda}}}^{\times}/\mu_{E_{\lambda}}$.
\end{cor}

By Proposition \ref{prop3.23.5}, we can assume (and we do) that $k(s)$ is the perfection of a finitely generated extension of $\mathbb{F}_p$, in which case the algebraic closure of $\mathbb{F}_p$ in $s$ is a finite field. Let $\F$ be an object of $M(X,1,\Ow_E)$. We prove the existence of a $d$-cycle $\mathcal{E}(\F)$ with the required properties, and with coefficients in the smaller abelian group
$$
A=\mathbb{Q} \otimes_{\Z} \Hom_{\mathrm{cont}}(W(\overline{s}/s), \Ow_E[p^{-1}]^{\times}/\mu_E),$$
where $\Ow_E[p^{-1}]^{\times}/\mu_E$ is endowed with the discrete topology. We note that for any element $\lambda$ of $\mathbb{L}$ not dividing $p$, the homomorphism
$$
A \rightarrow \mathbb{Q} \otimes_{\Z} \Hom_{\mathrm{cont}}(W(\overline{s}/s), \Ow_{E_{\lambda}}^{\times}/\mu_{E_{\lambda}}) \simeq \mathbb{Q} \otimes_{\Z} \Hom_{\mathrm{cont}}(\Gal(\overline{s}/s), \Ow_{E_{\lambda}}^{\times}/\mu_{E_{\lambda}}),
$$
is injective. Let $\F$ be an object of $M(X,1,\Ow_E)$. For any two elements $\lambda,\lambda'$ of $\mathbb{L}$ not dividing $p$, there exists a closed conical subset $C$ of $T^*X$ of pure dimension $d$ on which both $\F_{\lambda}$ and $\F_{\lambda'}$ are micro-supported. Let us consider a commutative diagram
\begin{center}
 \begin{tikzpicture}[scale=1]
\node (B) at (1,1) {$u$};
\node (A) at (2,0) {$S$};
\node (C) at (4,0) {$s$};
\node (D) at (2,1) {$U$};
\node (E) at (4,1) {$X$};

\path[->,font=\scriptsize]
(B) edge  (D)
(A) edge  (C)
(D) edge node[left]{$f$} (A)
(E) edge (C)
(D) edge node[above]{$j$} (E);
\end{tikzpicture} 
\end{center}
where $j$ is \'etale, $u$ is a closed point of $U$, $S$ is a smooth curve over $s$, and $df^{-1}(C_{|U \setminus \{ u \}})$ is contained in the $0$-section of $T^*S$. Let $\psi : \mathbb{F}_p \rightarrow E^{\times}$ be a non trivial homomorphism and let $\pi$ be a uniformizer on $S_{(f(u))}$. Then $\Art_{\pi,\psi } u^* R \Phi_{X \times_S S_{(f(u))}/ S_{(f(u))}} j^* \F$ is in $M(u, 1, E)$, cf. \ref{3.20} and Proposition \ref{propvareps}. By Proposition \ref{invar2} and by \ref{idento}, we obtain that the $\varepsilon$-factor of $(R \Phi_{X \times_S S_{(f(u))}/ S_{(f(u))}} j^* \F_{\lambda})_{\overline{u}}$ is in $A$ and coincides with the $\varepsilon$-factor of $(R \Phi_{X \times_S S_{(f(u))}/ S_{(f(u))}} j^* \F_{\lambda'})_{\overline{u}}$. Thus the intersection multiplicity $(\mathcal{E}(\F_{\lambda}),df)_{T^*U,u}$ is in $A$ and coincides with $(\mathcal{E}(\F_{\lambda'}),df)_{T^*U,u}$. This implies that $\mathcal{E}(\F_{\lambda})$ and $\mathcal{E}(\F_{\lambda'})$ have coefficients in $A$ and coincide.

\bibliographystyle{amsalpha}

\end{document}